\documentclass[a4paper,11pt,fleqn]{article}
\usepackage{amsmath,amssymb,amsthm,graphicx,subfigure,float,caption,epstopdf,tabularx,color, bm,amsfonts,epic}
\usepackage[top=1in, bottom=1in, left=1.25in, right=1.25in]{geometry}
\usepackage{appendix}
\usepackage{multirow}
\usepackage{diagbox}
\allowdisplaybreaks
\newtheorem{thm}{Theorem}[section]
\newtheorem{lem}{Lemma}[section]

\newtheorem{rmk}{Remark}[section]

\newtheorem*{prf}{Proof}
\numberwithin{equation}{section}

\usepackage{indentfirst}
\graphicspath{{Fig}}

\begin{document}
\title{A linearly implicit structure-preserving scheme for the Camassa-Holm equation based on multiple scalar auxiliary variables approach}
\author{Chaolong Jiang$^1$, Yuezheng Gong$^2$,\ Wenjun Cai$^3$,\ Yushun Wang$^3$\footnote{Correspondence author. Email:
wangyushun@njnu.edu.cn.}\\
{\small $^1$ School of Statistics and Mathematics, }\\
{\small Yunnan University of Finance and Economics, Kunming 650221, P.R. China}\\
{\small $^2$ College of Science,}\\
{\small Nanjing University of Aeronautics and Astronautics, Nanjing 210016, P.R. China}\\
{\small $^3$ Jiangsu Key Laboratory for Numerical Simulation of Large Scale Complex Systems,}\\
{\small School of Mathematical Sciences,  Nanjing Normal University,}\\
{\small  Nanjing 210023, P.R. China}\\
}
\date{}
\maketitle

\begin{abstract}
In this paper, we present a linearly implicit energy-preserving scheme for the Camassa-Holm equation by using the multiple scalar auxiliary variables approach, which is first developed to construct efficient and robust energy stable schemes for gradient systems. The Camassa-Holm equation is first reformulated into an
equivalent system by utilizing the multiple scalar auxiliary variables approach, which inherits
a modified energy. Then, the system is discretized in space aided by the standard Fourier pseudo-spectral method and a semi-discrete system is obtained, which is proven to preserve a semi-discrete modified energy. Subsequently, the linearized Crank-Nicolson method is applied for the resulting semi-discrete system to arrive at a fully discrete scheme. The main feature of the new
scheme is to form a linear system with a constant coefficient matrix at each time step and produce numerical solutions along which the modified energy is precisely conserved, as is the case with the analytical solution. Several numerical results are addressed to confirm accuracy and efficiency of the proposed scheme. \\[2ex]
\end{abstract}

\section{Introduction}
In this paper, we consider the Camassa-Holm (CH)  equation \cite{CH93,CHH94}
\begin{align}\label{ch:eq:1.1}
\left\lbrace
  \begin{aligned}
&u_t-u_{xxt}+3uu_x-2u_xu_{xx}-uu_{xxx}=0,\ a<x<b,\ 0<t\le T,\\
&u(x,0)=u_0(x),\ a\le x\le b,\\
&u(x+L,t)=u(x,t),\ a\le x\le b,\ 0\le t\le T,
\end{aligned}\right.
  \end{align}
where $t$ is time, $x$ is the spatial coordinate, $L=b-a$, and $u(x,t)$ represents the water's
free surface in non-dimensional variables. The CH equation models the unidirectional
propagation of shallow water waves over a flat bottom and is completely integrable \cite{CH93,Constantin01}. Thus, it has
infinitely many conservation laws. The first three are
\begin{align}\label{ch:eq:neq1.0}
&\frac{d}{dt}\mathcal{M}=0,\ \mathcal{M}=\int_{a}^b udx,\\\label{ch:eq:neq1.1}
&\frac{d}{dt}\mathcal{I}=0,\ \mathcal{I}=\int_{a}^b(u^2+u_x^2)dx,\\\label{ch:eq:neq1.2}
&\frac{d}{dt}\mathcal{H}=0,\ \mathcal{H}=-\frac{1}{2}\int_{a}^b\Big(u^3+uu_x^2\Big)dx,
\end{align}
where $\mathcal{M}$, $\mathcal{I}$ and $\mathcal{H}$ are the mass, momentum and energy of the CH equation \eqref{ch:eq:1.1}, respectively. The aim of this paper is concerned with the numerical methods that preserve the energy.

Because the energy is the most important first integral of the CH equation, designing of energy-preserving methods attracts a lot of interest. In Ref. \cite{MY09}, Matsuo et al. presented an energy-conserving Galerkin scheme for the CH equation. Further analysis for the energy-preserving $H^1$-Galerkin scheme was investigated in Ref. \cite{MM12}. Later on, Cohen and Raynaud \cite{CR11} derived a new energy-preserving scheme by the discrete gradient approach. Recently, Gong and Wang \cite{GW16b} proposed an energy-preserving wavelet collocation scheme for the CH equation \eqref{ch:eq:1.1}. However, such energy-preserving schemes are fully implicit that typically need iterative solvers at each time step. This quickly becomes a computationally expensive procedure. To address this drawback and maintain the desired energy-preserving property, Eidnes et al. \cite{ELS19} constructed two linearly implicit energy-preserving schemes for the CH equation \eqref{ch:eq:1.1} using the Kahan's method and the polarised discrete gradient methods, respectively. In Ref. \cite{JWG19}, we proposed a novel linearly implicit energy-preserving scheme for the CH equation \eqref{ch:eq:1.1} using the invariant energy quadratization (IEQ) approach \cite{GZYW18,YZW17,ZYGW17}. At each time step, the linearly implicit schemes only require to solve a linear system, which leads to considerably lower costs than the implicit one \cite{DO11}. However, these schemes leads to a linear system with complicated variable coefficients at each time step that may be difficult or expensive to solve.

More recently, inspired by the scalar auxiliary variable (SAV) approach \cite{SXY17,SXY18}, Cai et al. developed a linearly implicit energy-conserving scheme for the sine-Gordon equation \cite{CJWSjcp2019}. The resulting scheme leads to a linear system with constant coefficients that is easy to implement. The purpose of this paper is to apply the idea of the SAV approach to develop an efficient and energy-preserving scheme for the CH equation \eqref{ch:eq:1.1}. However, 
the classical SAV approach can not be directly applied to develop energy-preserving schemes for the CH equation. Actually, following the classical SAV approach, we introduce the auxiliary variable, as follows:
\begin{align}\label{ch:eq:new:1.4}
q=\sqrt{\int_{a}^b\Big(u^3+uu_x^2\Big)dx+C_0},
\end{align}
where $C_0$ is a constant large enough to make $q$ well-defined. The energy is then rewritten as
\begin{align}\label{ch:eq:1.3}
 \mathcal{H}=-\frac{1}{2}q^2+\frac{1}{2}C_0.
 \end{align}
 According to the energy variational, the CH equation \eqref{ch:eq:1.1} can be reformulated into an equivalent system, as follows:
\begin{align}\label{ch:eq:1.4}
\left\lbrace
  \begin{aligned}
  &\partial_t u=\mathcal{D}\Bigg(-\frac{3u^2+u_x^2}{2\sqrt{(u^3+uu_x^2,1)+C_0}}q+\frac{\partial_x\big(2uu_x\big)}{2\sqrt{(u^3+uu_x^2,1)+C_0}}q\Bigg),\\
  &\frac{d}{dt}q=0,\\
  &u(x,0)=u_0(x),\ q(0)=\sqrt{\int_{a}^b\Big(u_0(x)^3+u_0(x)\partial_xu_0(x)^2\Big)dx+C_0},\\
  &u(x+L,t)=u(x,t),
  \end{aligned}\right.\ \ 
  \end{align}
  where $\mathcal{D}=(1-\partial_{xx})^{-1}\partial_x$ is a skew-adjoint operator. {However, the above reformulated system has two main drawbacks for the further development of efficiently energy-preserving schemes: (i) according to the second equation of \eqref{ch:eq:1.4}, $q$ reduces a constant,  which fails to contribute to the numerical scheme; (ii) based on the conventional SAV discretization where $q$ is treated implicitly in time but other terms are treated explicitly, we obtain a fully explicit scheme, which may require a strict restriction on the grid ratio.} 
To meet these challenges, we first split the energy \eqref{ch:eq:neq1.2} as three parts, where two parts are bounded from below and the rest is quadratic. Then, we utilize the multiple scalar auxiliary variables (MSAV) approach \cite{CS18} to transform the original system into
an equivalent form, which inherits a modified energy. Subsequently, a novel linearly implicit energy-preserving scheme is proposed by applying the linearly implicit structure-preserving method in time and the standard Fourier pseudo-spectral method in space, respectively, for the equivalent system. We show that the proposed scheme can exactly preserve the discrete modified energy and mass, respectively, and only require to solve a linear system with a constant coefficient matrix at each time step that can be solved by FFT solvers efficiently. The MSAV approach is more recently {proposed by Cheng and Shen in Ref. \cite{CS18}} to deal with free energies with multiple disparate terms in the phase-field vesicle membrane and leads to robust energy stable schemes which enjoy the same computational advantages as the classical SAV approach. To the best of our knowledge, there is no result concerning the MSAV approach for the energy-conserving system. Taking the CH equation \eqref{ch:eq:1.1} for example, we first explore the feasibility of the MSAV approach and then devise a linearly implicit energy-preserving scheme. 

The outline of this paper is
organized as follows. In Section \ref{Sec:Ch:2}, based on the MSAV approach, the CH equation \eqref{ch:eq:1.1} is reformulated into an equivalent form.
A semi-discrete system, which inherits the semi-discrete modified energy, is presented in Section \ref{Sec:Ch:3}. In Section \ref{Sec:Ch:4},
we will concentrate on the construction for the linearly implicit energy-preserving scheme. Several numerical experiments are reported in Section \ref{Sec:Ch:5}. We draw some conclusions in Section \ref{Sec:Ch:6}.
\section{Model reformulation using the MSAV approach}\label{Sec:Ch:2}
In this section, we first reformulate the CH equation into an equivalent form with a
quadratic energy functional using the idea of the MSAV approach. The resulting reformulation provides an elegant platform for developing linearly implicit energy-preserving schemes.

The energy functional \eqref{ch:eq:neq1.2} can be split as the following three parts
  \begin{align}\label{ch:eq:2.1}
  \mathcal{H}
  &=-\frac{1}{2}\int_{a}^b(u+\frac{1}{2})^2(u^2+u_x^2)dx+\frac{1}{2}\int_{a}^bu^2(u^2+u_x^2)dx+\frac{1}{8}\int_{a}^b(u^2+u_x^2)dx\nonumber\\
  &:=-\frac{1}{2}\int_{a}^bg(u,u_x)dx+\frac{1}{2}\int_{a}^bh(u,u_x)dx+\frac{1}{8}\int_{a}^b(u^2+u_x^2)dx.
  \end{align}
Subsequently, following the idea of the MSAV approach, we introduce two scalar auxiliary variables, as follows:
{\begin{align*}
 q_1=\sqrt{(g(u,u_x),1)+C_1},\ q_2=\sqrt{(h(u,u_x),1)+C_2},
 \end{align*}}
 where $(v,w)$ is the inner product defined by $(v,w)=\int_{a}^b vw dx$, and $C_1$ and $C_2$ are two constants large enough to make $q_1$ and $q_2$ well-defined. Eq. \eqref{ch:eq:2.1} can then be rewritten as
 \begin{align}\label{ch:eq:2.2}
 \mathcal{H}=\frac{1}{8}\int_{a}^b(u^2+u_x^2)dx-\frac{1}{2}q_1^2+\frac{1}{2}q_2^2+\frac{1}{2}C_1-\frac{1}{2}C_2.
 \end{align}
 According to the energy variational, the system \eqref{ch:eq:1.1} can be reformulated into the following equivalent form
 \begin{align}\label{ch:eq:2.3}
\left\lbrace
  \begin{aligned}
  &\partial_t u=\mathcal{D}\Bigg(-\frac{\Big(\frac{\partial g}{\partial u}(u,u_x)-\partial_x\frac{\partial g}{\partial u_x}(u,u_x)\Big)}{2\sqrt{(g(u,u_x),1)+C_1}}q_1\\
  &~~~~~~~~~~~~+\frac{\Big(\frac{\partial h}{\partial u}(u,u_x)-\partial_x\frac{\partial h}{\partial u_x}(u,u_x)\Big)}{2\sqrt{(h(u,u_x),1)+C_2}}q_2+\frac{1}{4}(u-u_{xx})\Bigg),\\
  &\frac{d}{dt}q_1=\Bigg(\frac{1}{2\sqrt{(g(u,u_x),1)+C_1}}\Big(\frac{\partial g}{\partial u}(u,u_x)-\partial_x\frac{\partial g}{\partial u_x}(u,u_x)\Big),u_t\Bigg),\\
  &\frac{d}{dt}q_2=\Bigg(\frac{1}{2\sqrt{(h(u,u_x),1)+C_2}}\Big(\frac{\partial h}{\partial u}(u,u_x)-\partial_x\frac{\partial h}{\partial u_x}(u,u_x)\Big),u_t\Bigg),\\
  &u(x,0)=u_0(x),\ q_1(0)=\sqrt{(g(u_0(x),\partial_xu_0(x)),1)+C_1},\\
&q_2(0)=\sqrt{(h(u_0(x),\partial_xu_0(x)),1)+C_2},\\
&u(x+L,t)=u(x,t),
  \end{aligned}\right.\ \ 
  \end{align}
where
\begin{align*}
&\frac{\partial g}{\partial u}(u,u_x)=2(u+\frac{1}{2})(2u^2+u_x^2+\frac{1}{2}u),\ \frac{\partial g}{\partial u_x}=2u_x(u+\frac{1}{2})^2,\\
&\frac{\partial h}{\partial u}(u,u_x)=4u^3+2uu_x^2,\ \frac{\partial h}{\partial u_x}=2u_xu^2.
\end{align*}
\begin{thm} The system \eqref{ch:eq:2.3} possesses the following modified energy.
\begin{align*}
\frac{d}{dt}\mathcal{H}=0,\ \mathcal{H}=\frac{1}{8}\int_{a}^b(u^2+u_x^2)dx-\frac{1}{2}q_1^2+\frac{1}{2}q_2^2+\frac{1}{2}C_1-\frac{1}{2}C_2.
\end{align*}
\end{thm}
   \begin{prf}\rm We can deduce from \eqref{ch:eq:2.3} that
  \begin{align*}
   \frac{d}{dt} \mathcal{H}&=\frac{1}{4}(u-u_{xx},u_t)-q_1\frac{d}{dt}q_1+q_2\frac{d}{dt}q_2-\frac{1}{4}u_{x}u_t|_{a}^b\\
  &=\Bigg(-\frac{\Big(\frac{\partial g}{\partial u}(u,u_x)-\partial_x\frac{\partial g}{\partial u_x}(u,u_x)\Big)}{2\sqrt{(g(u,u_x),1)+C_1}}q_1
  +\frac{\Big(\frac{\partial h}{\partial u}(u,u_x)-\partial_x\frac{\partial h}{\partial u_x}(u,u_x)\Big)}{2\sqrt{(h(u,u_x),1)+C_2}}q_2\nonumber\\
  &~~~~~+\frac{1}{4}(u-u_{xx}),u_t\Bigg)-\frac{1}{4}u_{x}u_t|_{a}^b\\
  &=0,
   \end{align*}
   where 
   the last equality follows from the first equality of \eqref{ch:eq:2.3}, the periodic boundary condition and the skew-adjoint property of $\mathcal{D}$. \qed
   \end{prf}
\begin{rmk} We should note that the splitting strategy used in \eqref{ch:eq:2.1} is not unique. The comparisons between splitting strategies will be the subject of future investigations.
\end{rmk}

\section{Structure-preserving spatial semi-discretization}\label{Sec:Ch:3}
In this section, the standard Fourier pseudo-spectral method is employed to approximate spatial derivatives of the system \eqref{ch:eq:2.3} and we prove that
the resulting semi-discrete system can exactly preserve the semi-discrete modified energy.

Choose the mesh size $h=L/N$ with $N$ an even positive integer, and denote the grid points by $x_{j}=jh$ for $j=0,1,2,\cdots,N$; let $U_{j}$ be the numerical approximation of  $u(x_j,t)$ for $j=0,1,\cdots,N$ and ${U}=(U_{0},U_{1},\cdots,U_{N-1})^{T}$ be the solution vector space, and define discrete inner product as
\begin{align*}
&\langle {U},{ V}\rangle_{h}=h\sum_{j=0}^{N-1}U_{j}V_{j},\ \langle {U},{\bm 1}\rangle_{h}=h\sum_{j=0}^{N-1}U_{j},
\end{align*}
where {${\bm 1}=(1,1,\cdots,1)^T\in\mathbb{R}^N$.}

 Let
\begin{align*}
&S_{N}=\text{span}\{g_{j}(x),\ 0\leq j\leq N-1\}
\end{align*}
be the interpolation space, where $g_{j}(x)$ is trigonometric polynomials of degree $N/2$
given by
\begin{align*}
  &g_{j}(x)=\frac{1}{N}\sum_{l=-N/2}^{N/2}\frac{1}{a_{l}}e^{\text{i}l\mu (x-x_{j})},
\end{align*}
with $a_{l}=\left \{
 \aligned
 &1,\ |l|<\frac{N}{2},\\
 &2,\ |l|=\frac{N}{2},
 \endaligned
 \right.,$ and $\mu=\frac{2\pi}{b-a}$.
We define the interpolation operator $I_{N}: C(\Omega)\to S_{N}$, as follows:
\begin{align*}
I_{N}u(x,t)=\sum_{j=0}^{N-1}u_{j}(t)g_{j}(x),
\end{align*}
where $u_{j}(t)=u(x_{j},t)$. 
Taking the derivative with respect to $x$, and then evaluating the resulting expression at the collocation points $x_{j}$, we have
\begin{align*}
\frac{\partial^{s} I_{N}u(x_{j},t)}{\partial x^{s}}
&=\sum_{j_1=0}^{N-1}u_{j_1}\frac{d^{s}g_{j_1}(x_{j})}{dx^{s}}=[{D}_{s}{u}]_{j},
\end{align*}
where $u=(u_0,u_1,\cdots,u_{N-1})^T$ and ${D}_{s}$ is an $N\times N$ matrix, with elements given by
\begin{align*}
({ D}_{s})_{j_1,j}=\frac{d^{s}g_{j}(x_{j_1})}{dx^{s}}.
\end{align*}
In particular, the first and second order differential matrices can be obtained explicitly \cite{CQ01}
\begin{align*}
&({D}_{1})_{j,l}=
 \left \{
 \aligned
 &\frac{1}{2}\mu (-1)^{j+l} \cot(\mu \frac{x_{j}-x_{l}}{2}),\ &j\neq l,\\
 &0,\quad \quad \quad \quad \quad \quad \quad \ \ \ ~&j=l,
 \endaligned
 \right.
\end{align*}
\begin{align*}
 &({D}_{2})_{j,l}=
 \left \{
 \aligned
 &\frac{1}{2}\mu^{2} (-1)^{j+l+1}\csc^{2}(\mu \frac{x_{j}-x_{l}}{2}),\ &j\neq l,\\
 &-\mu^{2}\frac{N^{2}+2}{12},\quad \quad \quad \quad ~ &j=l.
 \endaligned
 \right.
 \end{align*}
 \begin{rmk}\label{ch:rmk:3.1} It should be remarked that, for matrices $D_1$ and $D_2$, the following result holds \cite{GCW14,ST06}
 \begin{align*}
&{D}_{1}={F}_{N}^{H}\Lambda_1{F}_{N},\ \Lambda_1=\text{\rm i}\mu\text{\rm diag}\Big[0,1,\cdots,\frac{N}{2}-1,0,1-\frac{N}{2},\cdots,-1\Big],\\
&{D}_{2}={F}_{N}^{H}\Lambda_2{F}_{N},\ \Lambda_2=-\mu^2\text{\rm diag}\Big[0^2,1^2,\cdots,\big(\frac{N}{2}\big)^2,\big(-\frac{N}{2}+1\big)^2,\cdots,(-2)^2,(-1)^2\Big],
\end{align*}
 where  ${F}_{N}$ is the discrete Fourier transform matrix with elements
$\big({F}_{N}\big)_{j,k}=\frac{1}{\sqrt{N}}e^{-\text{\rm i}jk\frac{2\pi}{N}},\ 0\le j,k\le N-1,$ ${F}_{N}^{H}$ is the conjugate transpose matrix of ${F}_{N}$.
 \end{rmk}
 Applying the standard Fourier pseudo-spectral method to the system \eqref{ch:eq:2.3} in space and we have
 \begin{align}\label{ch:eq:3.4}
\left\lbrace
  \begin{aligned}
  &\frac{d}{dt} {U}={D}\Bigg(-\frac{\Big(g_1({U},{D}_1{U})-{D}_1 g_2({U},{D}_1{ U})\Big)}{2\sqrt{\langle g({U},{D}_1{U}),{\bm 1}\rangle_h+C_1}}q_1\\
  &~~~~~~~~~+\frac{\Big(h_1({U},{D}_1{U})-{D}_1h_2({U},{D}_1{U})\Big)}{2\sqrt{\langle h({U},{D}_1{U}),{\bm 1}\rangle_h+C_2}}q_2+\frac{1}{4}\big({U}-{D}_2{U}\big)\Bigg),\\
  &\frac{d}{dt} q_1=\Bigg\langle \frac{\Big(g_1({U},{D}_1{U})-{D}_1 g_2({U},{D}_1{ U})\Big)}{2\sqrt{\langle g({U},{D}_1{U}),{\bm 1}\rangle_h+C_1}},\frac{d}{dt} { U}\Bigg\rangle_h,\\
   &\frac{d}{dt} q_2=\Bigg\langle \frac{\Big(h_1({U},{D}_1{U})-{D}_1h_2({U},{D}_1{ U})\Big)}{2\sqrt{\langle h({U},{D}_1{U}),{\bm 1}\rangle_h+C_2}},\frac{d}{dt} {U}\Bigg\rangle_h,\\
  \end{aligned}\right.\ \ 
  \end{align}
  where ${D}=({I}-{D}_2)^{-1}{D}_1$, $g_1=\frac{\partial g}{\partial u},\ g_2=\frac{\partial g}{\partial u_x}$, $h_1=\frac{\partial h}{\partial u}$, and $h_2=\frac{\partial h}{\partial u_x}$. 

 \begin{thm} The semi-discrete system \eqref{ch:eq:3.4} admits the following semi-discrete modified energy
\begin{align*}
\frac{d}{dt} E_h=0,\ E_h=\frac{1}{8}\langle {U}-{D}_2{U},{U}\rangle_h-\frac{1}{2}q_1^2+\frac{1}{2}q_2^2+\frac{1}{2}C_1-\frac{1}{2}C_2.
\end{align*}
\end{thm}
   \begin{prf}\rm It follows from the semi-discrete system \eqref{ch:eq:3.4} that
   \begin{align*}
   \frac{d}{dt} E_h&=\frac{1}{4}\langle {U}-{D}_2{U},\frac{d}{dt}{U}\rangle_h-q_1\frac{d}{dt} q_1+q_2\frac{d}{dt} q_2\\
   &=\Big\langle -\frac{\Big(g_1({U},{D}_1{U})-{D}_1 g_2({U},{D}_1{U})\Big)}{2\sqrt{\langle g({U},{D}_1{U}),{\bm 1}\rangle_h+C_1}}q_1+\frac{\Big(h_1({U},{D}_1{U})-{D}_1h_2({U},{D}_1{U})\Big)}{2\sqrt{\langle h({U},{D}_1{U}),{\bm 1}\rangle_h+C_2}}q_2\\
   &~~~~+\frac{1}{4}({U}-{D}_2{U}),\frac{d}{dt} {U}\Big\rangle_h\\
   &=0,
   \end{align*} where
the last equality follows from the first equality \eqref{ch:eq:3.4} and the skew-symmetry of ${D}$.
   \qed
   \end{prf}

 \section{Construction of the linearly implicit energy-preserving scheme }\label{Sec:Ch:4}
 In this section, we present a linearly implicit energy-preserving scheme by utilizing the linearized Crank-Nicolson method to the semi-discrete system \eqref{ch:eq:3.4} in time.

Choose $\tau=T/M$ be the time step with $M$ a positive integer, and denote $t_{n}=n\tau$ for $n=0,1,2\cdots,M$; {let $U_{j}^n,\ Q_1^n$ and $Q_2^n$ be the numerical approximations of  $u(x_j,t_n),\ q_1(t_n)$ and $q_2(t_n)$, respectively, for $j=0,1,\cdots,N$ and $n=0,1,2,\cdots,M$; denote $U^n=(U_0^n,U_1^n,\cdots,U_{N-1}^n)^T$ }as the solution
vector at $t=t_n$ and define
\begin{align*}
&\delta_{t} {U}_j^{n}
=\frac{{U}_j^{n+1}-{U}_j^{n}}{\tau},\ {U}_j^{n+\frac{1}{2}}=\frac{{U}_j^{n+1}+{U}_j^{n}}{2},\ \hat{U}_j^{n+\frac{1}{2}}=\frac{3{U}_j^{n}-{U}_j^{n-1}}{2}, 0\le j\le N-1.
\end{align*}

 Applying the linearized Crank-Nicolson method to the semi-discrete system \eqref{ch:eq:3.4} in time, and we obtain a fully discretized scheme, as follows:
 \begin{align}\label{ch:eq:4.1}
\left\lbrace
  \begin{aligned}
  &\delta_t {U}^n={D}\Bigg(-\frac{\Big(g_1(\hat{U}^{n+\frac{1}{2}},{D}_1\hat{U}^{n+\frac{1}{2}})-{D}_1 g_2(\hat{U}^{n+\frac{1}{2}},{D}_1\hat{U}^{n+\frac{1}{2}})\Big)}{2\sqrt{\langle g(\hat{U}^{n+\frac{1}{2}},{D}_1\hat{ U}^{n+\frac{1}{2}}),{\bm 1}\rangle_h+C_1}}Q_1^{n+\frac{1}{2}}\\
  &~~~~~~~~~+\frac{\Big(h_1(\hat{U}^{n+\frac{1}{2}},{D}_1\hat{U}^{n+\frac{1}{2}})-{D}_1h_2(\hat{U}^{n+\frac{1}{2}},{D}_1\hat{U}^{n+\frac{1}{2}})\Big)}{2\sqrt{\langle h(\hat{U}^{n+\frac{1}{2}},{D}_1\hat{ U}^{n+\frac{1}{2}}),{\bm  1}\rangle_h+C_2}}Q_2^{n+\frac{1}{2}}\\
  &~~~~~~~~~+\frac{1}{4}\big({U}^{n+\frac{1}{2}}-{D}_2{U}^{n+\frac{1}{2}}\big)\Bigg),\\
  &\delta_t Q_1^n=\Bigg\langle \frac{\Big(g_1(\hat{U}^{n+\frac{1}{2}},{D}_1\hat{U}^{n+\frac{1}{2}})-{D}_1 g_2(\hat{U}^{n+\frac{1}{2}},{D}_1\hat{U}^{n+\frac{1}{2}})\Big)}{2\sqrt{\langle g(\hat{U}^{n+\frac{1}{2}},{D}_1\hat{ U}^{n+\frac{1}{2}}),{\bm 1}\rangle_h+C_1}},\delta_t {U}^n\Bigg\rangle_h,\\
   &\delta_t Q_2^n=\Bigg\langle \frac{\Big(h_1(\hat{U}^{n+\frac{1}{2}},{D}_1\hat{U}^{n+\frac{1}{2}})-{D}_1h_2(\hat{U}^{n+\frac{1}{2}},{D}_1\hat{U}^{n+\frac{1}{2}})\Big)}{2\sqrt{\langle h(\hat{U}^{n+\frac{1}{2}},{D}_1\hat{ U}^{n+\frac{1}{2}}),{\bm 1}\rangle_h+C_2}},\delta_t {U}^n\Bigg\rangle_h,\\
  \end{aligned}\right.\ \ 
  \end{align}
 for $n=1,\cdots,M-1$. { Since the scheme \eqref{ch:eq:4.1} is three-level, we obtain ${U}^1,\ Q_1^1$ and $Q_2^1$ by
 \begin{align}\label{ch:eq:4.2}
\left\lbrace
  \begin{aligned}
  &\delta_t {U}^0={D}\Bigg(-\frac{\Big(g_1({U}^{0},{D}_1{U}^{0})-{D}_1 g_2({U}^{0},{D}_1{U}^{0})\Big)}{2\sqrt{\langle g({U}^{0},{D}_1{ U}^{0}),{\bm 1}\rangle_h+C_1}}Q_1^{\frac{1}{2}}\\
  &~~~~~~~~~+\frac{\Big(h_1({U}^{0},{D}_1{U}^{0})-{D}_1h_2({U}^{0},{D}_1{U}^{0})\Big)}{2\sqrt{\langle h({U}^{0},{D}_1{ U}^{0}),{\bm  1}\rangle_h+C_2}}Q_2^{\frac{1}{2}}+\frac{1}{4}\big({U}^{\frac{1}{2}}-{D}_2{U}^{\frac{1}{2}}\big)\Bigg),\\
  &\delta_t Q_1^0=\Bigg\langle \frac{\Big(g_1({U}^{0},{D}_1{U}^{0})-{D}_1 g_2({U}^{0},{D}_1{U}^{0})\Big)}{2\sqrt{\langle g({U}^{0},{D}_1{ U}^{0}),{\bm 1}\rangle_h+C_1}},\delta_t {U}^0\Bigg\rangle_h,\\
   &\delta_t Q_2^0=\Bigg\langle \frac{\Big(h_1({U}^{0},{D}_1{U}^{0})-{D}_1h_2({U}^{0},{D}_1{U}^{0})\Big)}{2\sqrt{\langle h({U}^{0},{D}_1{ U}^{0}),{\bm 1}\rangle_h+C_2}},\delta_t {U}^0\Bigg\rangle_h.\\
  \end{aligned}\right.\ \ 
  \end{align}}
 The initial and boundary conditions in \eqref{ch:eq:2.3} are discretized as
\begin{align*}
&U_j^0=u_0(x_j),\ Q_1^0=\sqrt{\langle g(U^0,D_1U^0),{\bm 1}\rangle_h+C_1},\ Q_2^0=\sqrt{\langle h(U^0,D_1U^0),{\bm 1}\rangle_h+C_2},\\
& U_{j\pm N}^n=U_j^n,\ j=0,1,2,\cdots,N.
\end{align*}
%

Then, we show that the proposed scheme \eqref{ch:eq:4.1}-\eqref{ch:eq:4.2} can exactly preserve the discrete energy and  mass, respectively.
\begin{thm} The proposed scheme \eqref{ch:eq:4.1}-\eqref{ch:eq:4.2} satisfies the following discrete modified energy
\begin{align}\label{ch:eq-4.4}
E_h^{n+1}=E_h^n,\ E_h^n=\frac{1}{8}\langle {U}^n-{D}_2{U}^n,{U}^n\rangle_h-\frac{1}{2}(Q_1^n)^2+\frac{1}{2}(Q_2^n)^2+\frac{1}{2}C_1-\frac{1}{2}C_2,
\end{align}
for $n=0,1,\cdots,M-1.$
\end{thm}
   \begin{prf}\rm It is readily to obtain from \eqref{ch:eq:4.1} that
   \begin{align}\label{ch:eq:new4.3}
   \delta_t E_h^n&=\frac{1}{4}\langle {U}^{n+\frac{1}{2}}-{D}_2{ U}^{n+\frac{1}{2}},\delta_t{U}^n \rangle_h-Q_1^{n+\frac{1}{2}}\delta_tQ_1^n+Q_2^{n+\frac{1}{2}}\delta_t Q_2^n\nonumber\\
   &=\Big\langle -\frac{\Big(g_1(\hat{U}^{n+\frac{1}{2}},{D}_1\hat{U}^{n+\frac{1}{2}})-{D}_1 g_2(\hat{U}^{n+\frac{1}{2}},{D}_1\hat{ U}^{n+\frac{1}{2}})\Big)}{2\sqrt{\langle g(\hat{U}^{n+\frac{1}{2}},{D}_1\hat{U}^{n+\frac{1}{2}}),{\bm 1}\rangle_h+C_1}}Q_1^{n+\frac{1}{2}}\nonumber\\
   &~~+\frac{\Big(h_1(\hat{U}^{n+\frac{1}{2}},{D}_1\hat{U}^{n+\frac{1}{2}})-{D}_1h_2(\hat{U}^{n+\frac{1}{2}},{D}_1\hat{ U}^{n+\frac{1}{2}})\Big)}{2\sqrt{\langle h(\hat{U}^{n+\frac{1}{2}},{D}_1\hat{U}^{n+\frac{1}{2}}),{\bm 1}\rangle_h+C_2}}Q_2^{n+\frac{1}{2}}\nonumber\\
   &~~+\frac{1}{4}({U}^{n+\frac{1}{2}}-{D}_2{U}^{n+\frac{1}{2}}),\delta_t {U}^n\Big\rangle_h\nonumber\\
   &=0,
   \end{align}
   which further implies
   \begin{align*}
   E_h^{n+1}=E_h^n,\ n=1,2,\cdots,M-1,
   \end{align*}
    where the last equality of \eqref{ch:eq:new4.3} follows from the first equality of \eqref{ch:eq:4.1} and the skew-symmetry of ${D}$.
   An argument similar to \eqref{ch:eq:4.2} used in \eqref{ch:eq:new4.3} shows that
   \begin{align*}
   E_h^{1}=E_h^0.
   \end{align*}
   This completes the proof.
   \qed
   \end{prf}
{\begin{thm} The scheme \eqref{ch:eq:4.1}-\eqref{ch:eq:4.2}  possesses the discrete mass
\begin{align*}
M_h^{n+1}=M_h^n,\ M_h^n=\langle U^n,{\bm 1}\rangle_h,\ n=0,1,2,\cdots,M-1.
\end{align*}
\end{thm}}
\begin{prf}\rm {According to Remark \ref{ch:rmk:3.1}, for $j=0,1,2,\cdots,N-1$, we have
\begin{align*}
(F_N{\bm 1})_j=\frac{1}{\sqrt{N}}\sum_{k=0}^{N-1}e^{-\text{\rm i}jk\frac{2\pi}{N}}=
 \left \{
 \aligned
 &\sqrt{N},\ &\  j=0,\\
 &0,\ &1\le j\le N-1,
 \endaligned
 \right.
\end{align*}
and the elements of matrix $(I-\Lambda_2)^{-1}\Lambda_1$ in the first row are all zero.}

{Taking the discrete inner product of \eqref{ch:eq:4.1} with ${\bm 1}$, we then obtain
\begin{align*}
\langle\delta_t {U}^n,{\bm 1}\rangle_h&=\Bigg\langle{D}\Bigg(-\frac{\Big(g_1(\hat{U}^{n+\frac{1}{2}},{D}_1\hat{U}^{n+\frac{1}{2}})-{D}_1 g_2(\hat{U}^{n+\frac{1}{2}},{D}_1\hat{U}^{n+\frac{1}{2}})\Big)}{2\sqrt{\langle g(\hat{U}^{n+\frac{1}{2}},{D}_1\hat{ U}^{n+\frac{1}{2}}),{\bm 1}\rangle_h+C_1}}Q_1^{n+\frac{1}{2}}\\
  &~~~~~~~~~+\frac{\Big(h_1(\hat{U}^{n+\frac{1}{2}},{D}_1\hat{U}^{n+\frac{1}{2}})-{D}_1h_2(\hat{U}^{n+\frac{1}{2}},{D}_1\hat{U}^{n+\frac{1}{2}})\Big)}{2\sqrt{\langle h(\hat{U}^{n+\frac{1}{2}},{D}_1\hat{ U}^{n+\frac{1}{2}}),{\bm  1}\rangle_h+C_2}}Q_2^{n+\frac{1}{2}}\\
  &~~~~~~~~~+\frac{1}{4}\big({U}^{n+\frac{1}{2}}-{D}_2{U}^{n+\frac{1}{2}}\big)\Bigg),{\bm 1}\Bigg\rangle_h\\
  &=\Bigg\langle (I-\Lambda_2)^{-1}\Lambda_1F_N\Bigg(-\frac{\Big(g_1(\hat{U}^{n+\frac{1}{2}},{D}_1\hat{U}^{n+\frac{1}{2}})-{D}_1 g_2(\hat{U}^{n+\frac{1}{2}},{D}_1\hat{U}^{n+\frac{1}{2}})\Big)}{2\sqrt{\langle g(\hat{U}^{n+\frac{1}{2}},{D}_1\hat{ U}^{n+\frac{1}{2}}),{\bm 1}\rangle_h+C_1}}Q_1^{n+\frac{1}{2}}\\
  &~~~~~~~~~+\frac{\Big(h_1(\hat{U}^{n+\frac{1}{2}},{D}_1\hat{U}^{n+\frac{1}{2}})-{D}_1h_2(\hat{U}^{n+\frac{1}{2}},{D}_1\hat{U}^{n+\frac{1}{2}})\Big)}{2\sqrt{\langle h(\hat{U}^{n+\frac{1}{2}},{D}_1\hat{ U}^{n+\frac{1}{2}}),{\bm  1}\rangle_h+C_2}}Q_2^{n+\frac{1}{2}}\\
  &~~~~~~~~~+\frac{1}{4}\big({U}^{n+\frac{1}{2}}-{D}_2{U}^{n+\frac{1}{2}}\big)\Bigg),F_N{\bm 1}\Bigg\rangle_h=0,
\end{align*}
which further shows
\begin{align*}
M_h^{n+1}=M_h^n,\ n=1,2,\cdots,M-1.
\end{align*}
By an argument similar to \eqref{ch:eq:4.1} used as above, we obtain
\begin{align*}
M_h^{1}=M_h^0.
\end{align*}
This completes proof. \qed}
\end{prf}

  Besides its energy-preserving property, a most remarkable thing about the above scheme is that it can be solved efficiently. Let
  \begin{small}
  \begin{align*}
  &{G}_1=\frac{1}{2\sqrt{\langle g(\hat{U}^{n+\frac{1}{2}},{D}_1\hat{U}^{n+\frac{1}{2}}),{1}\rangle_h+C_1}}\Big(g_1(\hat{ U}^{n+\frac{1}{2}},{D}_1\hat{U}^{n+\frac{1}{2}})-{D}_1 g_2(\hat{U}^{n+\frac{1}{2}},{D}_1\hat{U}^{n+\frac{1}{2}})\Big),\\
  &{G}_2=\frac{1}{2\sqrt{\langle h(\hat{U}^{n+\frac{1}{2}},{D}_1\hat{U}^{n+\frac{1}{2}}),{1}\rangle_h+C_2}}\Big(h_1(\hat{ U}^{n+\frac{1}{2}},{D}_1\hat{U}^{n+\frac{1}{2}})-{D}_1h_2(\hat{U}^{n+\frac{1}{2}},{D}_1\hat{U}^{n+\frac{1}{2}})\Big).
  \end{align*}
  \end{small}
   Eq. \eqref{ch:eq:4.1} can then rewritten as
  \begin{align}\label{ch:eq:4.3}
\left\lbrace
  \begin{aligned}
  &{U}^{n+\frac{1}{2}}={U}^n+\frac{\tau}{2}{D}\Bigg(-{G}_1Q_1^{n+\frac{1}{2}}+{G}_2Q_2^{n+\frac{1}{2}}+\frac{1}{4}\big({U}^{n+\frac{1}{2}}-{D}_2{ U}^{n+\frac{1}{2}}\big)\Bigg),\\
  &Q_1^{n+\frac{1}{2}}=Q_1^n+\Big\langle {G}_1,{U}^{n+\frac{1}{2}}-{U}^n\Big\rangle_h,\\
   &Q_2^{n+\frac{1}{2}}=Q_2^n+\Big\langle {G}_2,{U}^{n+\frac{1}{2}}-{U}^n\Big\rangle_h.\\
  \end{aligned}\right.\ \ 
  \end{align}
Next, by eliminating $Q_1^{n+\frac{1}{2}}$ and $Q_2^{n+\frac{1}{2}}$ from \eqref{ch:eq:4.3}, we have
  \begin{align}\label{ch:eq:4.4}
  \Big[{I}-\frac{\tau}{8}{D}_1&\Big]{U}^{n+\frac{1}{2}}=-\frac{\tau}{2}{D}{G}_1\langle {G}_1,{U}^{n+\frac{1}{2}}\rangle_h+\frac{\tau}{2}{D}{G}_2\langle { G}_2,{U}^{n+\frac{1}{2}}\rangle_h+{r}^n,
  \end{align}
  where
  \begin{align*}
  {r}^n&={U}^n-\frac{\tau}{2}{D}{G}_1Q_1^n+\frac{\tau}{2}{D}{G}_2Q_2^n+\frac{\tau}{2}{D}{G}_1\langle {G}_1,{U}^{n}\rangle_h-\frac{\tau}{2}{D}{G}_2\langle { G}_2,{U}^{n}\rangle_h.
  \end{align*}
Denote ${A}^{-1}=({I}-\frac{\tau}{8}{D}_1)^{-1}$ and
\begin{align*}
  &{\gamma}_1^n=-\frac{\tau}{2}{A}^{-1}{D}{G}_1,\
  {\gamma}_2^n=\frac{\tau}{2}{A}^{-1}{D}{G}_2, \ {b}^n={A}^{-1}{r}^n,
  \end{align*}
the above equation is equivalent to
  \begin{align}\label{ch:eq:4.5}
  {U}^{n+\frac{1}{2}}
  &={\gamma}_1^n\langle {G}_1,{U}^{n+\frac{1}{2}}\rangle_h+{\gamma}_2^n\langle {G}_2,{U}^{n+\frac{1}{2}}\rangle_h+{b}^n.
  \end{align}
  We take the inner product of \eqref{ch:eq:4.5} with ${G}_1$ and have
   \begin{align}\label{ch:eq:4.6}
  \langle {G}_1,{U}^{n+\frac{1}{2}}\rangle_h=\langle {G}_1,{\gamma}_1^n\rangle_h\langle {G}_1,{U}^{n+\frac{1}{2}}\rangle_h+\langle {G}_1,{ \gamma}_2^n\rangle_h\langle {G}_2,{U}^{n+\frac{1}{2}}\rangle_h+\langle {G}_1,{b}^n\rangle_h.
  \end{align}
Taking the inner product of \eqref{ch:eq:4.5} with ${G}_2$, we then obtain
   \begin{align}\label{ch:eq:4.7}
  \langle {G}_2,{U}^{n+\frac{1}{2}}\rangle_h=\langle {G}_2,{\gamma}_1^n\rangle_h\langle {G}_1,{U}^{n+\frac{1}{2}}\rangle_h+\langle {G}_2,{ \gamma}_2^n\rangle_h\langle {G}_2,{U}^{n+\frac{1}{2}}\rangle_h+\langle {G}_2,{b}^n\rangle_h.
  \end{align}
  Eqs. \eqref{ch:eq:4.6} and \eqref{ch:eq:4.7} form a $2\times 2$ linear system for the unknowns $(\langle {G}_1,{U}^{n+\frac{1}{2}}\rangle_h, \langle {G}_2,{U}^{n+\frac{1}{2}}\rangle_h)^T$.

Solving $(\langle {G}_1,{U}^{n+\frac{1}{2}}\rangle_h, \langle {G}_2,{U}^{n+\frac{1}{2}}\rangle_h)^T$ from the $2\times 2$ linear system \eqref{ch:eq:4.6} and \eqref{ch:eq:4.7} and ${U}^{n+\frac{1}{2}}$ is then updated from \eqref{ch:eq:4.5}. Subsequently, $Q_1^{n+\frac{1}{2}}$ and $Q_2^{n+\frac{1}{2}}$ are obtained from the second and third equality of \eqref{ch:eq:4.3}, respectively. Finally, we have ${U}^{n+1}=2{U}^{n+\frac{1}{2}}-{U}^{n}$, $Q_1^{n+1}=2Q_1^{n+\frac{1}{2}}-Q_1^{n}$ and $Q_2^{n+1}=2Q_2^{n+\frac{1}{2}}-Q_2^{n}$.
\begin{rmk} We should remark that, compared with the scheme obtained by the classical SAV approach, the proposed scheme need to solve an additional $2\times 2$ linear system, however, the main computational cost still comes from \eqref{ch:eq:4.4}. Thus, our scheme enjoys the
same computational advantages as the ones obtained by the classical SAV approach. In addition, in our computation, ${U}^{n+\frac{1}{2}}$ can be efficiently obtained from \eqref{ch:eq:4.5} by the FFT, when ones note Remark \ref{ch:rmk:3.1}.
\end{rmk}
 \begin{rmk} 
We should note that the energy \eqref{ch:eq:2.2} is equivalent to
the energy \eqref{ch:eq:neq1.2} in continuous sense, but not for the discrete sense. This indicates
that the scheme \eqref{ch:eq:4.1} cannot preserve the following discrete energy
\begin{align}\label{ch-eq:4.10}
H^n=-\frac{h}{2}\sum_{j=0}^{N-1}\Big((U_j^n)^3+U_j^n\cdot({D}_1{U}^n)_j^2\Big),\ 0\le n\le M.
\end{align}
\end{rmk}

\section{Numerical examples}\label{Sec:Ch:5}
In this section, we report the numerical performance, accuracy, CPU time and
invariants-preserving properties of the proposed scheme \eqref{ch:eq:4.1} (denoted by MSAV-LCNS). In addition, the following structure-preserving schemes are chosen for comparisons:
\begin{itemize}
\item IEQ-LCNS: the linearly implicit energy-preserving scheme given in Ref. \cite{JWG19};
\item  EPFPS: the energy-preserving Fourier pseudo-spectral scheme;
\item MSFPS: the multi-symplectic Fourier pseudo-spectral
scheme;
\item LICNS: the linear-implicit Crank-Nicolson scheme described in Ref. \cite{HGL19};
\item LILFS: the leap-frog scheme stated in Ref. \cite{HGL19}.
\end{itemize}
It is noted that EPFPS and MSFPS are obtained by using the Fourier pseudo-spectral method instead of the wavelet collocation method in Refs. \cite{GW16b,ZST11} , respectively. As a summary, a detailed table on the properties of each scheme has been given in Tab. \ref{Tab_CH:1}.

 In our computation, the FFT is also adopt as the solver of linear
systems given by MSAV-LCNS (see \eqref{ch:eq:4.4}), the standard fixed-point iteration is used for the fully
implicit schemes, and the Jacobi iteration method is employed for the linear systems given by IEQ-LCNS, LICNS and LILFS. Here, the iteration
will terminate if the infinity norm of the error between two adjacent iterative steps is less
than $10^{-14}$. All diagrams presented below refer to the numerical integration of the CH equation \eqref{ch:eq:1.1} with parameters $C_1=0$ and $C_2=0$, and are carried out via Matlab 7.0 with AMD A8-7100 and RAM
4GB. In order to quantify the numerical solution, we use the $l^2$- and $l^{\infty}$-norms of the error between the numerical solution $U_j^n$ and the exact solution
$u(x_j,t_n)$, respectively, as
\begin{align*}
e_{h,2}^2(t_n)=h\sum_{j=0}^{N-1}|U_{j}^n-u(x_j,t_n)|^2,\ e_{h,\infty}(t_n)=\max\limits_{0\le j\le N-1}|U_{j}^n-u(x_j,t_n)|,\ n\ge 0.
\end{align*}
\begin{table}[H]
\tabcolsep=6pt
\small
\renewcommand\arraystretch{1.1}
\centering
\caption{Comparison of properties of different numerical schemes}\label{Tab_CH:1}
\begin{tabular*}{\textwidth}[h]{@{\extracolsep{\fill}}c c c c c c c c c }\hline 
 \diagbox{Property}{Scheme}& MSAV-LCNS& IEQ-LCNS& EPFPS & MSFPS&LICNS&LILFS\\\hline
 Symplectic &No&No&No&Yes&No&No   \\ [1ex]  
 Energy conservation& Yes&Yes&Yes&No&Yes&Yes \\[1ex]
 Momentum conservation&No&No&No&No&Yes&Yes \\[1ex]
 Mass conservation&Yes&Yes&Yes&Yes&No&No \\[1ex]
 Fully implicit&No&No&Yes&Yes&No&No\\[1ex]
 Linearly implicit&Yes&Yes&No&No&Yes&Yes\\[1ex]
\hline
\end{tabular*}
\end{table}


\subsection{Smooth periodic solution}
{ The CH equation \eqref{ch:eq:1.1} admits smooth periodic traveling wave solutions
\begin{align*}
u(x,t)=\phi(x-ct),\ c\in\mathbb{R},
\end{align*}
when three parameters $m,M,z$ fulfill the relation $z<m<M<c$, where $m=\min\limits_{x\in\mathbb{R}}\phi,\ M=\max\limits_{x\in\mathbb{R}}\phi$ and $z=c-M-m$ (see Ref. \cite{lenellsjde05}). However, the solutions can only be given implicitly by
\begin{align}\label{ch:eq:5.1}
|x-x_0|=\int_{\phi_0}^{\phi}\frac{\sqrt{c-y}}{\sqrt{(M-y)(y-m)(y-z)}}dy,
\end{align}
where $\phi(x_0)=\phi_0$. In order to remove the singularities of the integral at $y=m$ and $y=M$, we transform \eqref{ch:eq:5.1} into
\begin{align}\label{ch:eq:5.2}
|x-x_0|=2\int_{\theta_0}^{\theta}\frac{\sqrt{A-\sin^2t}}{\sqrt{B+\sin^2 t}}dt
\end{align}
by the change of variables
\begin{align*}
\phi=m+(M-m)\sin^2\theta,
\end{align*}
where $A=\frac{c-m}{M-m}$ and $B=\frac{m-z}{M-m}$. The exact
solution is obtained by periodic extension of the initial data, which is constructed as follows:
\begin{itemize}
\item[Step 1:] Setting $x_0=0$ and $\theta_0=0$ and computing a grid on the interval $[0,\pi]$ with $\mathbb{N}$ equispaced nodes $\theta_j$ given by
\begin{align*}
\theta_j=\frac{\pi}{\mathbb{N}}j,\ j=0,1,2,\cdots,\mathbb{N}-1.
\end{align*}
\item[Step 2:] Computing $\phi_j=m+(M-m)\sin^2\theta_j,\ j=0,1,2,\cdots,\mathbb{N}-1$ and $x_j$ is obtained by using the Gaussian-Legendre quadrature formula for the integral $2\int_{0}^{\theta_j}\frac{\sqrt{A-\sin^2t}}{\sqrt{B+\sin^2 t}}dt,\ j=0,1,2,\cdots,\mathbb{N}-1$.
\item[Step 3:] Performing a spline interpolation for points $(x_j,\phi_j),\ j=0,1,2,\cdots,\mathbb{N}-1$ to obtain $\phi$ as a function of $x$.
\end{itemize}}

We take the bounded computational domain as the interval [$0,L$] with a periodic boundary condition and choose parameters $m=0.3, M=0.7$ and $c=1$, which gives rise to a smooth traveling wave with with period $L\approx 6.56$ (see Ref. \cite{KL05}).

\begin{figure}[H]
\centering\begin{minipage}[t]{70mm}
\includegraphics[width=70mm]{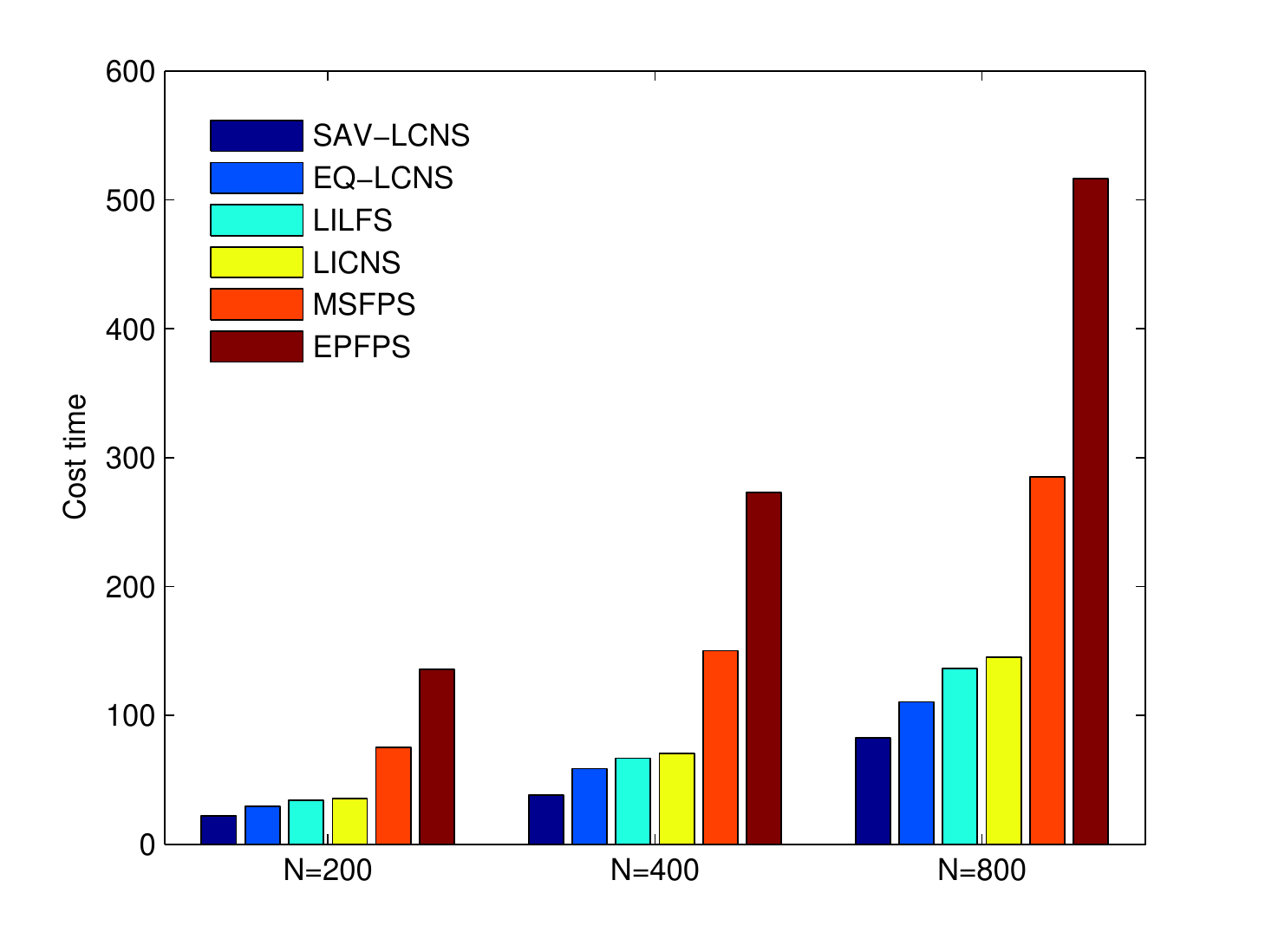}
\end{minipage}
\caption{The CPU times of the six schemes for the smooth solution under different mesh points till
$t=6.56$ with $\tau$=6.56e-04. }\label{Fig:ch:1}
\end{figure}

To test the temporal discretization errors of the different numerical schemes, we fix the Fourier node $32$ such that the spatial discretization errors are negligible. Tab. \ref{Tab:ch:1}, shows the temporal errors and convergence rates for different numerical schemes under different time steps at $t=6.56$.
Fig. \ref{Fig:ch:1} shows the CPU times of the six schemes for the smooth solution under different grid points till
$t=6.56$ with $\tau$=6.56e-04. From Tab. \ref{Tab:ch:1} and Fig. \ref{Fig:ch:1}, we can draw the following observations: (i) all schemes have second order accuracy in time errors; (ii) the error provided by MSAV-LCNS has the same order of magnitude as the ones provided by IEQ-LCNS and LICNS. (iii) the costs of EPFPS is
most expensive while the one of MSAV-LCNS is cheapest.

\begin{table}[H]
\tabcolsep=9pt
\footnotesize
\renewcommand\arraystretch{1.1}
\centering
\caption{{The numerical errors and convergence rates for different numerical schemes under different time steps at $t=6.56$.}}\label{Tab:ch:1}
\begin{tabular*}{\textwidth}[h]{@{\extracolsep{\fill}} c l l l l l}\hline
{Scheme\ \ } &{$\tau$} &{$e_{h,2}$} &{order}& {$e_{h,\infty}$}&{order}  \\     
\hline
 \multirow{4}{*}{MSAV-LCNS}  &{$\frac{L}{200}$}& {2.132e-03}&{-} &{1.485e-03} & {-}\\[1ex]
  {}  &{$\frac{L}{400}$}& {5.309e-04}&{2.01} &{3.717e-04} & {2.00}\\[1ex] 
   {}  &{$\frac{L}{800}$}& {1.327e-04}&{2.00} &{9.318e-05} &{2.00} \\[1ex]
    {}  &{$\frac{L}{1600}$}& {3.322e-05}&{2.00} &{2.334e-05} &{2.00} \\\hline
 \multirow{4}{*}{IEQ-LCNS}  &{$\frac{L}{200}$}& {1.050e-03}&{-} &{7.385e-04} & {-}\\[1ex]
  {}  &{$\frac{L}{400}$}& {2.600e-04}&{2.01} &{1.840e-04} & {2.01}\\ [1ex]
   {}  &{$\frac{L}{800}$}& {6.479e-05}&{2.00} &{4.598e-05} &{2.00} \\[1ex]
   {} &{$\frac{L}{1600}$} & {1.619e-05}&{2.00}&{1.150e-05} & {2.00}\\\hline
 \multirow{4}{*}{EPFPS}  &{$\frac{L}{200}$}& {4.005e-04}&{-} &{2.810e-04} & {-}\\[1ex]
  {}  &{$\frac{L}{400}$}& {1.002e-04}&{2.00} &{7.034e-05} & {2.00}\\[1ex] 
   {}  &{$\frac{L}{800}$}& {2.504e-05}&{2.00} &{1.759e-05} &{2.00} \\[1ex]
  {} &{$\frac{L}{1600}$} & {6.239e-06}&{2.00}&{4.396e-06} & {2.00}\\\hline
 \multirow{4}{*}{MSFPS}  &{$\frac{L}{200}$}& {4.106e-04}&{-} &{2.840e-04} & {-}\\[1ex]
  {}  &{$\frac{L}{400}$}& {1.027e-04}&{2.00} &{7.109e-05} & {2.00}\\[1ex] 
   {}  &{$\frac{L}{800}$}& {2.567e-05}&{2.00} &{1.778e-05} &{2.00} \\[1ex]
    {}  &{$\frac{L}{1600}$}& {6.397e-06}&{2.00} &{4.443e-06} &{2.00} \\\hline
    \multirow{4}{*}{LICNS}  &{$\frac{L}{200}$}& {8.884e-04}&{-} &{6.690e-04} & {-}\\[1ex]
  {}  &{$\frac{L}{400}$}& {2.213e-04}&{2.01} &{1.668e-04} & {2.00}\\[1ex] 
   {}  &{$\frac{L}{800}$}& {5.532e-05}&{2.00} &{4.171e-05} &{2.00} \\[1ex]
    {}  &{$\frac{L}{1600}$}& {1.385e-05}&{2.00} &{1.042e-05} &{2.00} \\\hline
    \multirow{4}{*}{LILFS}  &{$\frac{L}{200}$}& {4.447e-04}&{-} &{3.382e-04} & {-}\\[1ex]
  {}  &{$\frac{L}{400}$}& {1.111e-04}&{2.00} &{8.451e-05} & {2.01}\\[1ex] 
   {}  &{$\frac{L}{800}$}& {2.774e-05}&{2.00} &{2.110e-05} &{2.00} \\[1ex]
    {}  &{$\frac{L}{1600}$}& {6.903e-06}&{2.00} &{5.240e-06} &{2.00} \\\hline
\end{tabular*}
\end{table}


\begin{figure}[H]
\centering\begin{minipage}[t]{60mm}
\includegraphics[width=60mm]{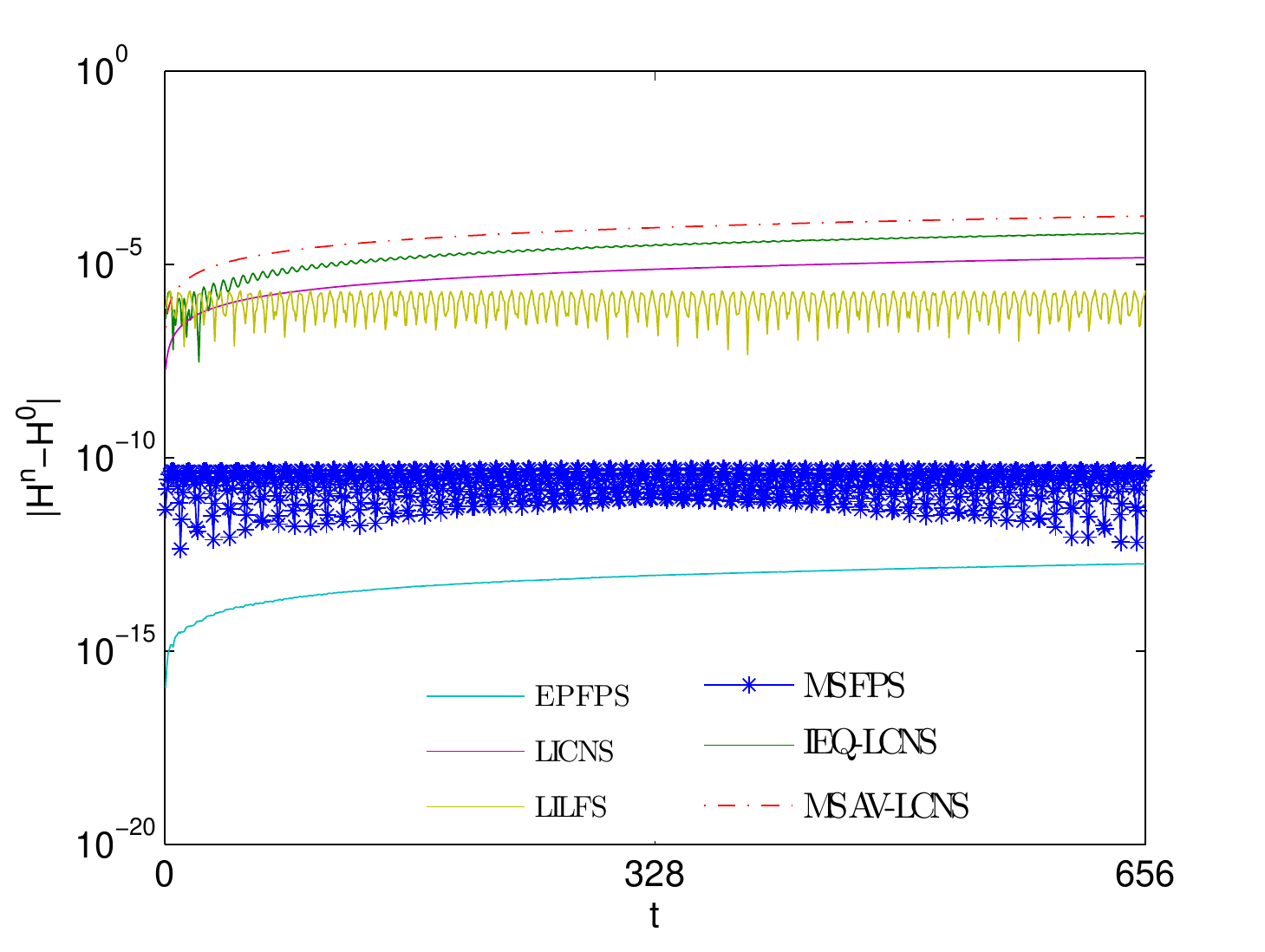}
\caption*{(a) Hamiltonian energy}
\end{minipage}\ \
\begin{minipage}[t]{60mm}
\includegraphics[width=60mm]{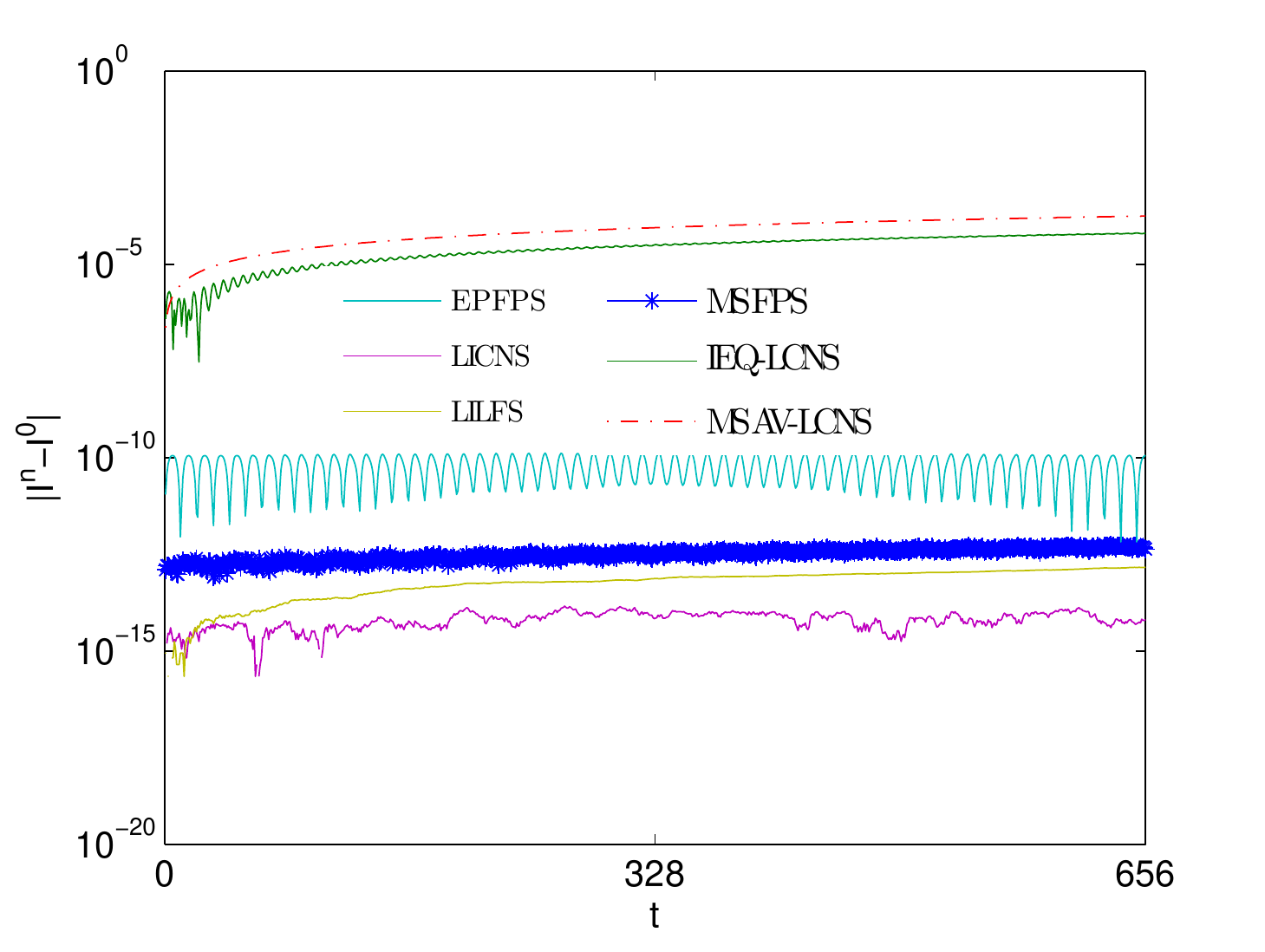}
\caption*{(b) Momentum}
\end{minipage}
\centering\begin{minipage}[t]{60mm}
\includegraphics[width=60mm]{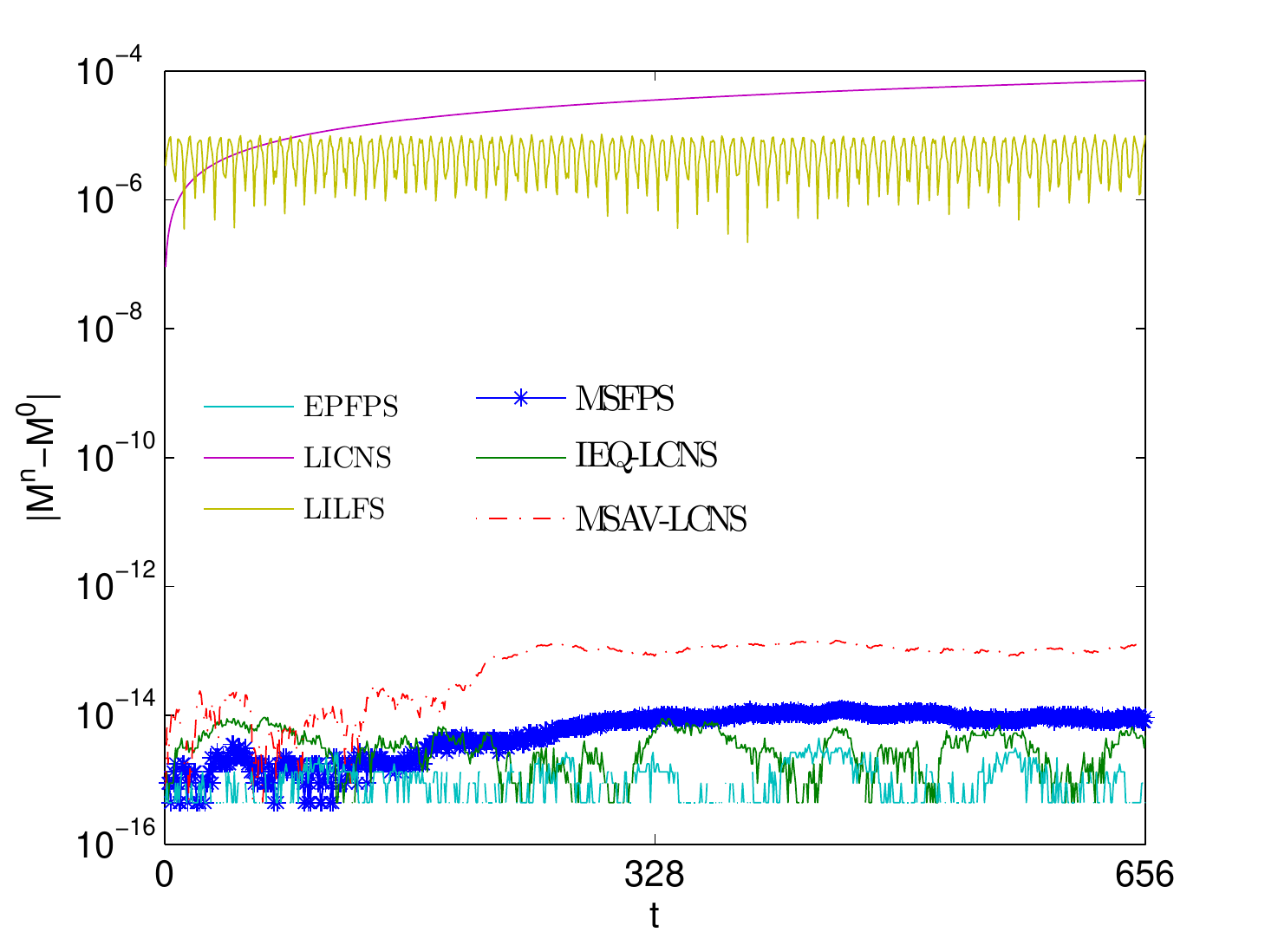}
\caption*{(c) Mass}
\end{minipage}\ \
\begin{minipage}[t]{60mm}
\includegraphics[width=60mm]{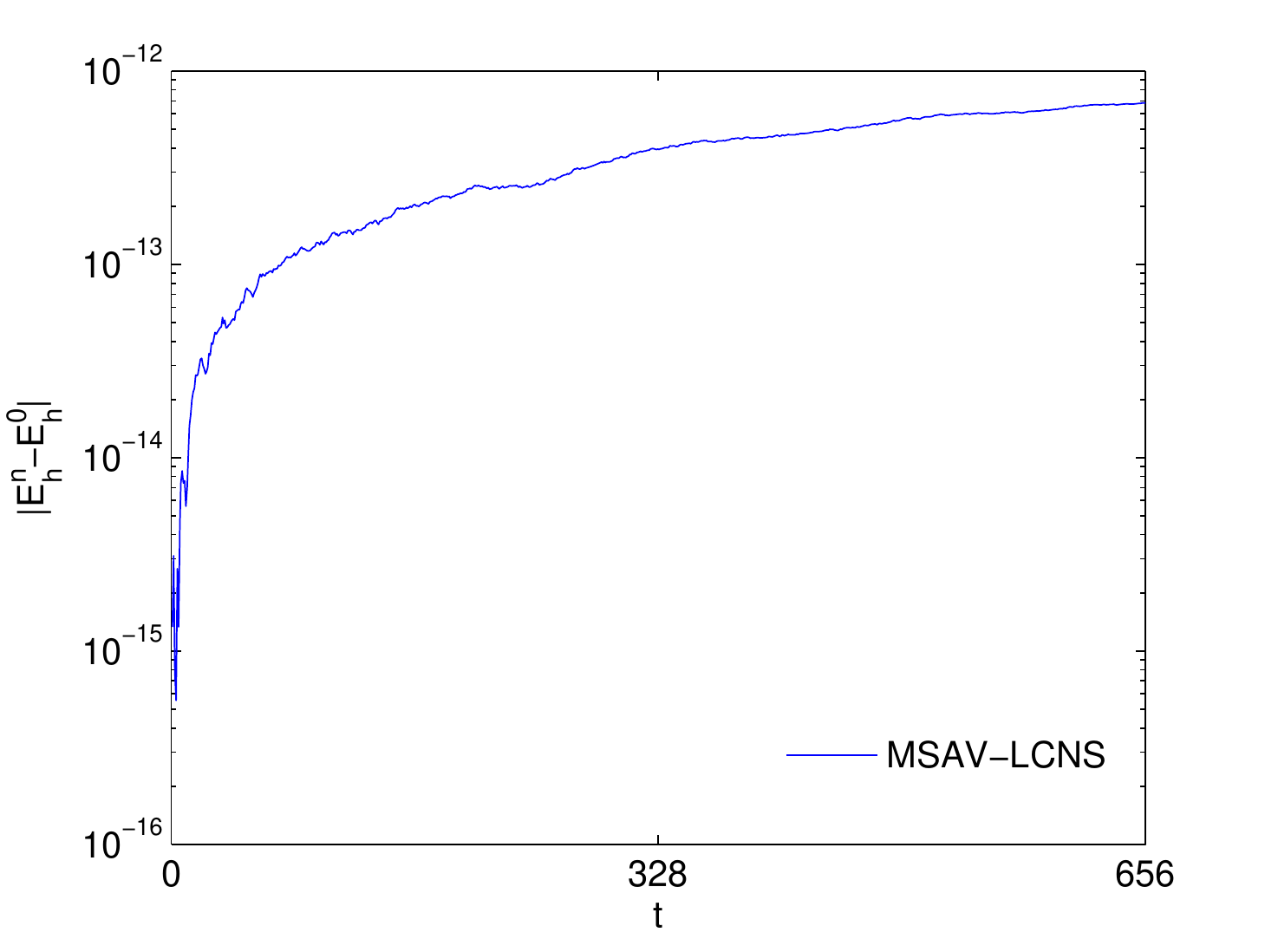}
\caption*{(d) Energy \eqref{ch:eq-4.4}}
\end{minipage}
\caption{The errors in invariants under $N$=32 and $\tau$=0.0082 over the time interval $t\in[0,656]$.}\label{Fig:ch:2}
\end{figure}
To further investigate the invariants-preservation of the proposed scheme. Fig. \ref{Fig:ch:2} shows the errors of the invariants under $N$=32 and $\tau$=0.0082 over the time interval $t\in[0,656]$. From Fig. \ref{Fig:ch:2} (a)-(c), we make the following observations: (i) EPFPS can exactly preserve the energy (see \eqref{ch-eq:4.10}) and the energy errors of others are remained around a small order of magnitude. (ii) LICNS and LILFS can exactly preserve the momentum and MSAV-LCNS, IEQ-LCNS, EPFPS and MSFPS can preserve the momentum approximately. (iii) MSAV-LCNS, IEQ-LCNS, EPFPS and MSFPS can preserve the mass to round-of errors while LICNS and LILFS admit large errors. From Fig. \ref{Fig:ch:2} (d), it is clearly demonstrated that the proposed scheme can exactly preserve the discrete modified energy. Similar observations on the errors of the invariants are made in the next three examples and we will omit these details for brevity. Here, we should note that the modified energy \eqref{ch:eq-4.4} and the energy \eqref{ch-eq:4.10} are an approximate version of the continue energy \eqref{ch:eq:neq1.2}, and the errors show the stability and capability for long-term computation of the numerical scheme.

\subsection{Two-peakon interaction}

\begin{figure}[H]
\centering
\centering\begin{minipage}[t]{60mm}
\includegraphics[width=65mm]{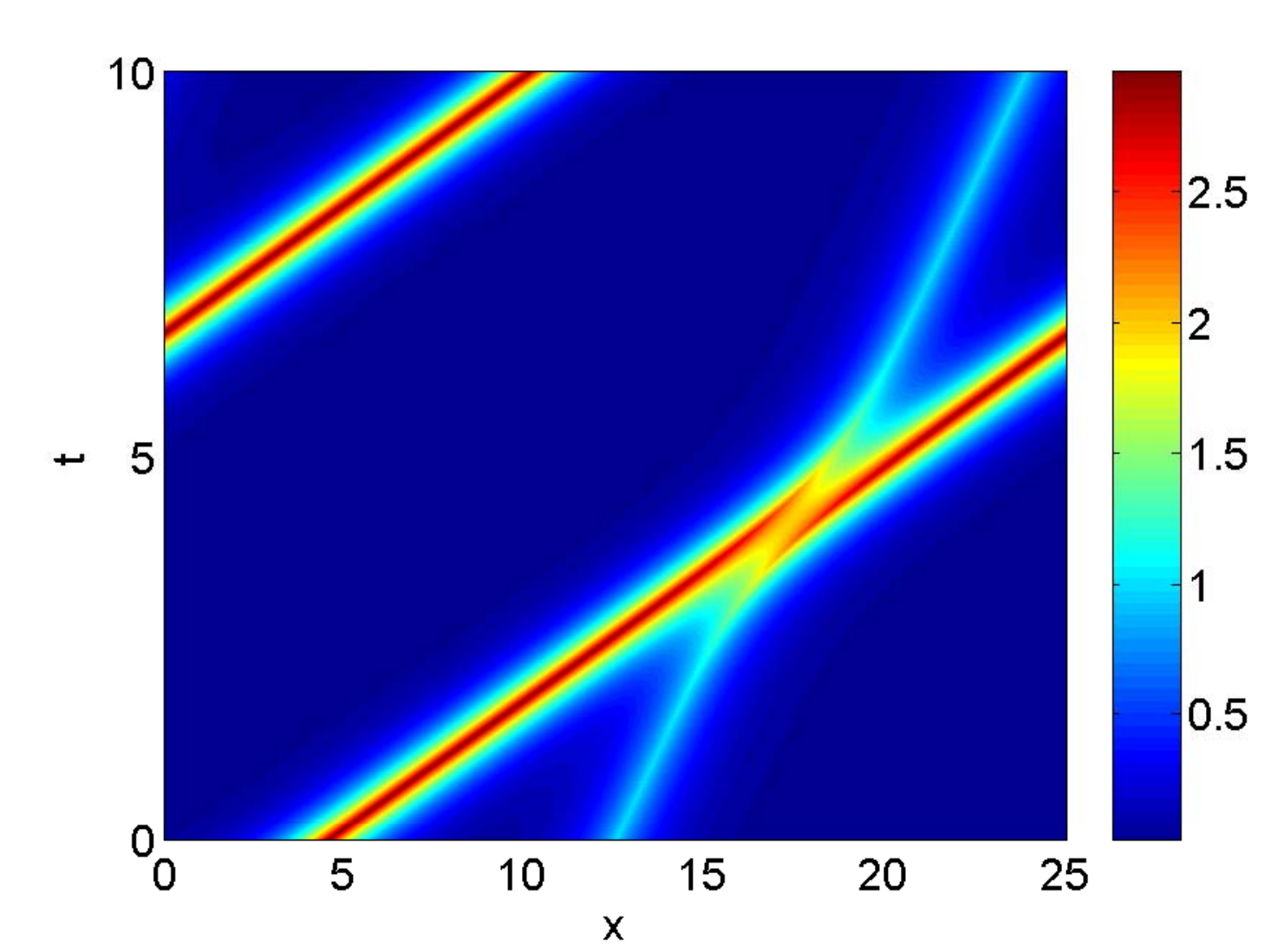}
\end{minipage}\ \
\caption{The two-peakon interaction of the CH equation \eqref{ch:eq:1.1} under $N$=1024 and $\tau=0.0001$.}\label{Fig:ch:3}
\end{figure}

\begin{figure}[H]
\centering\begin{minipage}[t]{60mm}
\includegraphics[width=65mm]{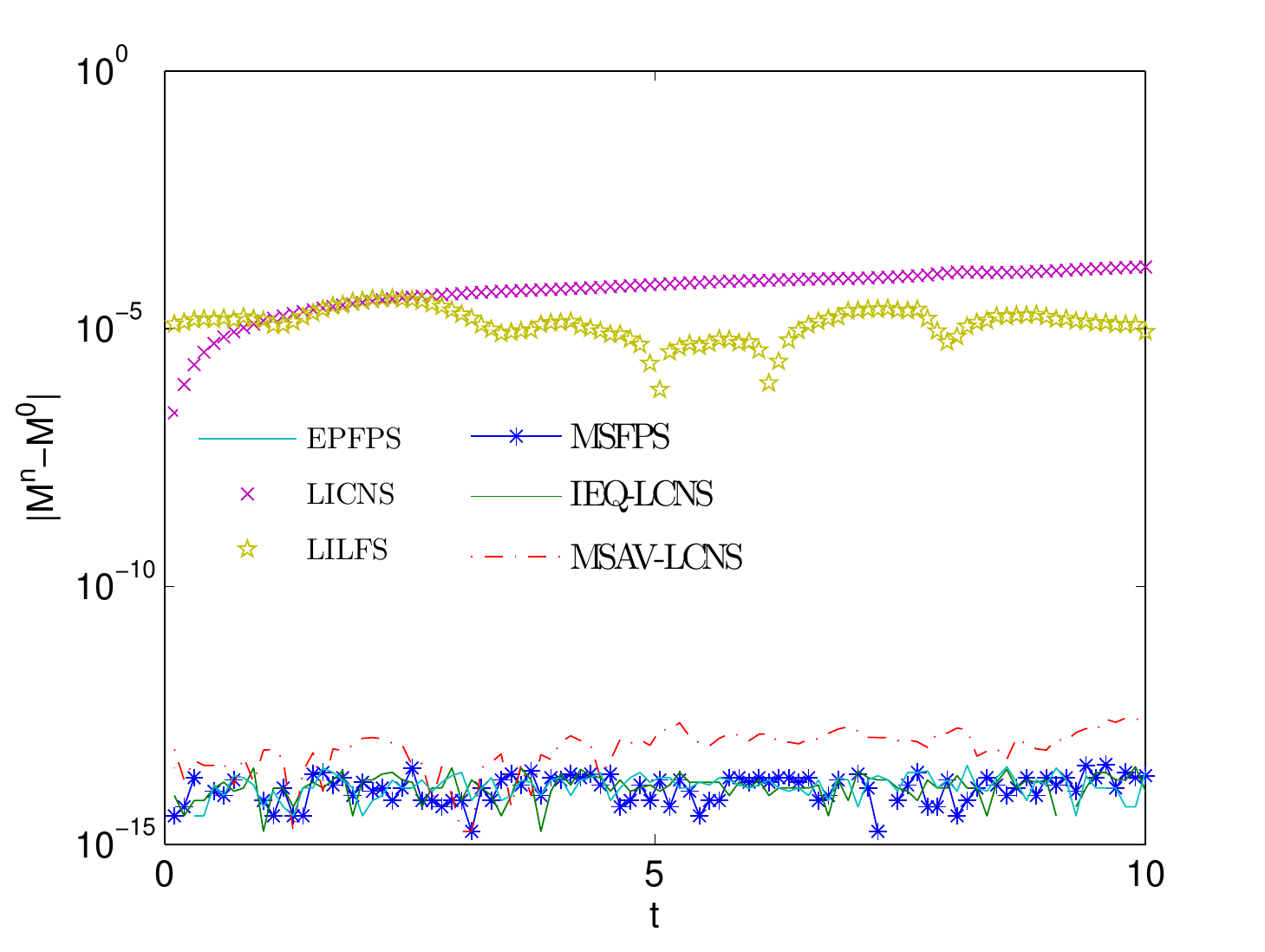}
\caption*{(c) Mass}
\end{minipage}\ \
\begin{minipage}[t]{60mm}
\includegraphics[width=65mm]{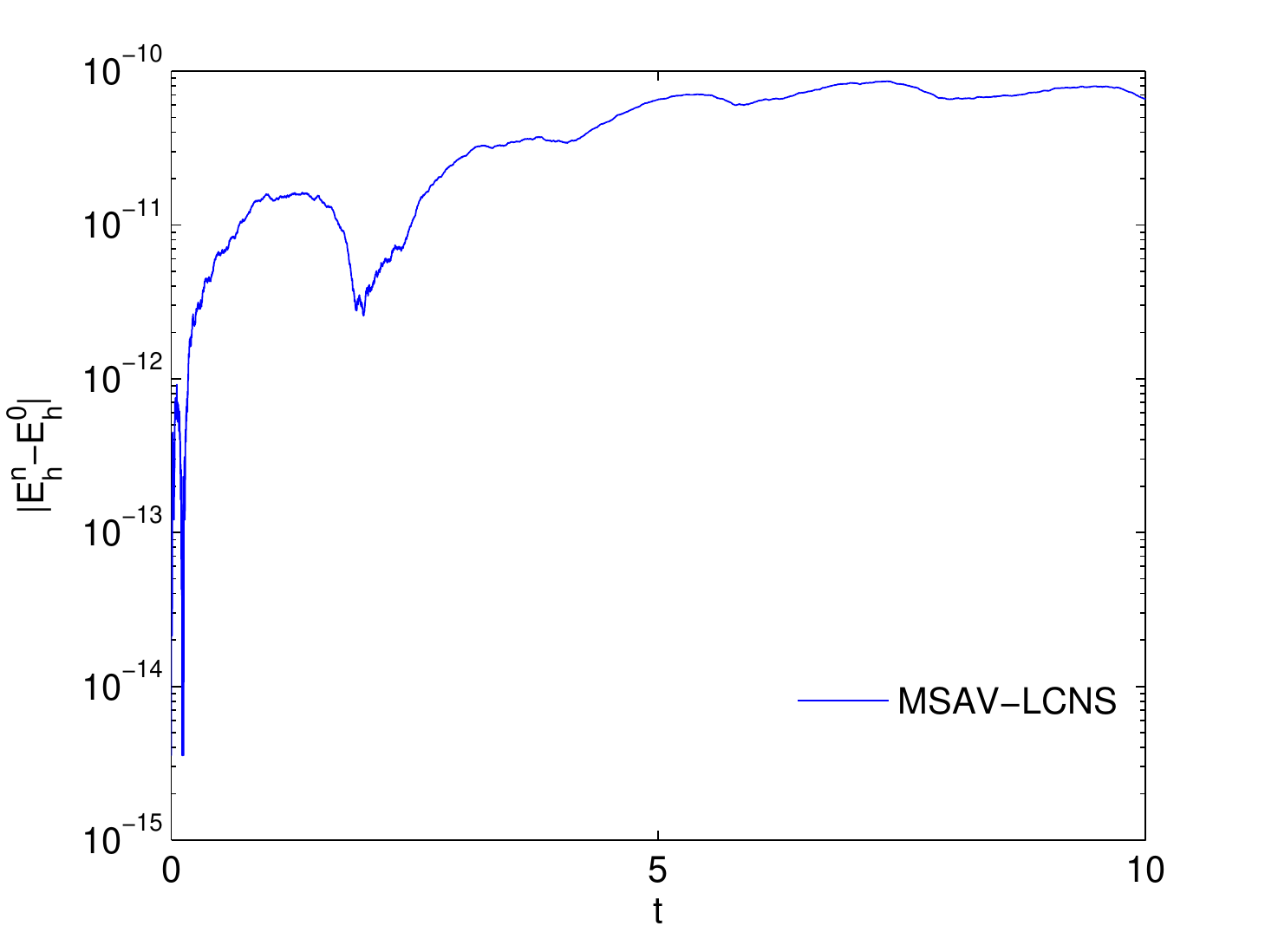}
\caption*{(d) Energy \eqref{ch:eq-4.4}}
\end{minipage}
\caption{The errors in invariants under $N$=1024 and $\tau$=0.0001 over the time interval $t\in[0,10]$.}\label{Fig:ch:4}
\end{figure}

We consider the two-peakon interaction of the CH equation \eqref{ch:eq:1.1} with the initial condition \cite{XS08}
\begin{align*}
u_0(x)=\phi_1(x)+\phi_2(x),\ 0\le x\le 25,
\end{align*}
where
\begin{align*}
\phi_i(x)=\left\lbrace
  \begin{aligned}
  &\frac{c_i}{\cosh(L/2)}\cosh(x-x_i),\ |x-x_i|\le L/2,\\
  &\frac{c_i}{\cosh(L/2)}\cosh(L-(x-x_i)),\ |x-x_i|> L/2,\
  \end{aligned} i=1,2.
  \right.\ \ 
  \end{align*}
The parameters are $c_1=3,c_2=1,x_1=-8,x_2=0,L=25$ and a periodic
boundary condition is considered. Fig. \ref{Fig:ch:3} shows the contour plot of the two-peakon interaction. We can
see clearly that the taller wave overtakes the shorter one and afterwards both waves retain their original shapes and velocities. The errors of invariants under $N$=1024 and $\tau$=0.0001 over the time interval $t\in[0,10]$ are plotted in Fig. \ref{Fig:ch:4}, which behaves similarly as that of Fig. \ref{Fig:ch:2}.

\begin{figure}[H]
\centering\begin{minipage}[t]{60mm}
\includegraphics[width=65mm]{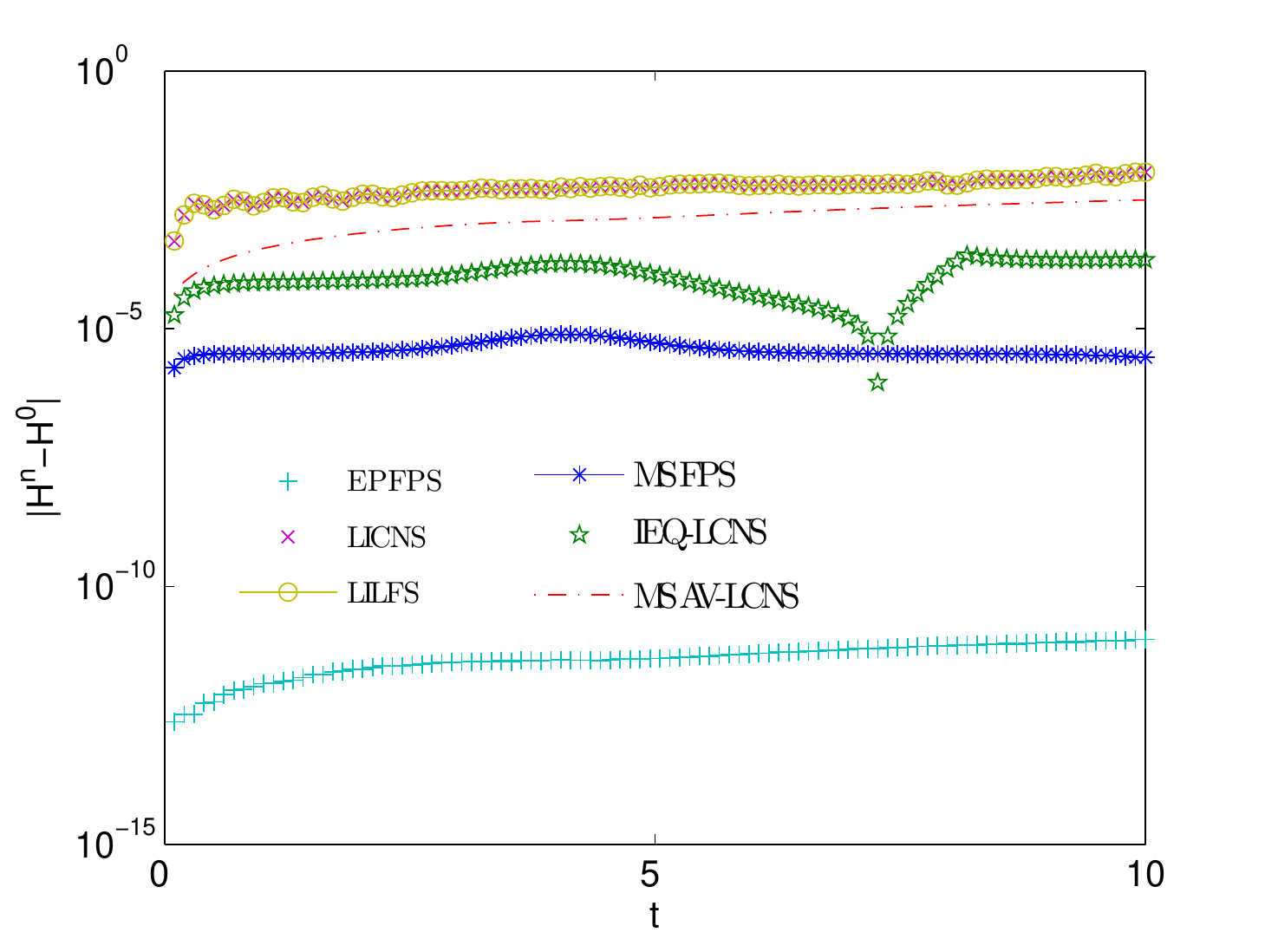}
\caption*{(a) Hamiltonian energy}
\end{minipage}\ \
\begin{minipage}[t]{60mm}
\includegraphics[width=65mm]{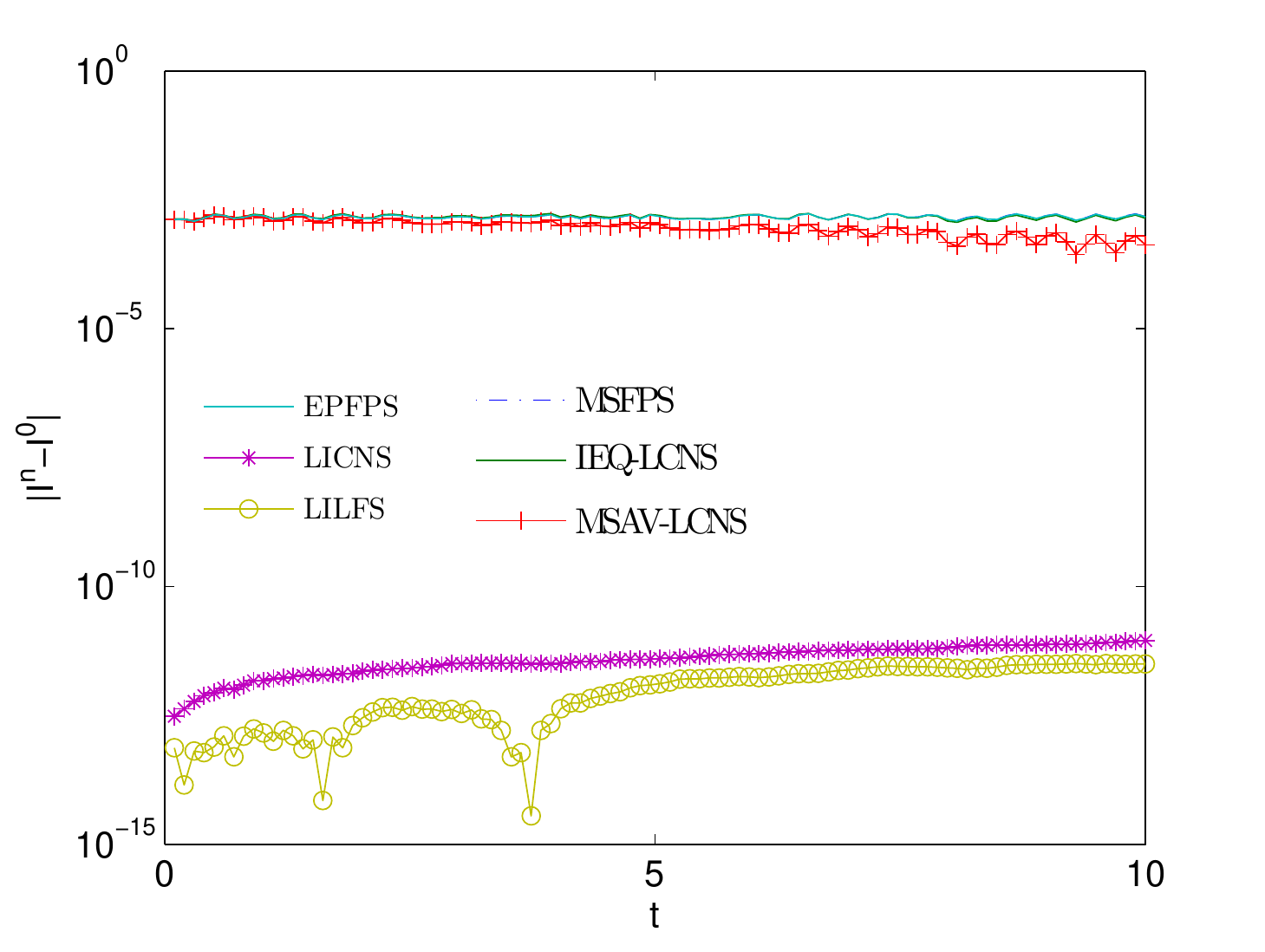}
\caption*{(b) Momentum}
\end{minipage}
\end{figure}

\subsection{Three-peakon interaction}
\begin{figure}[H]
\centering\begin{minipage}[t]{80mm}
\includegraphics[width=80mm]{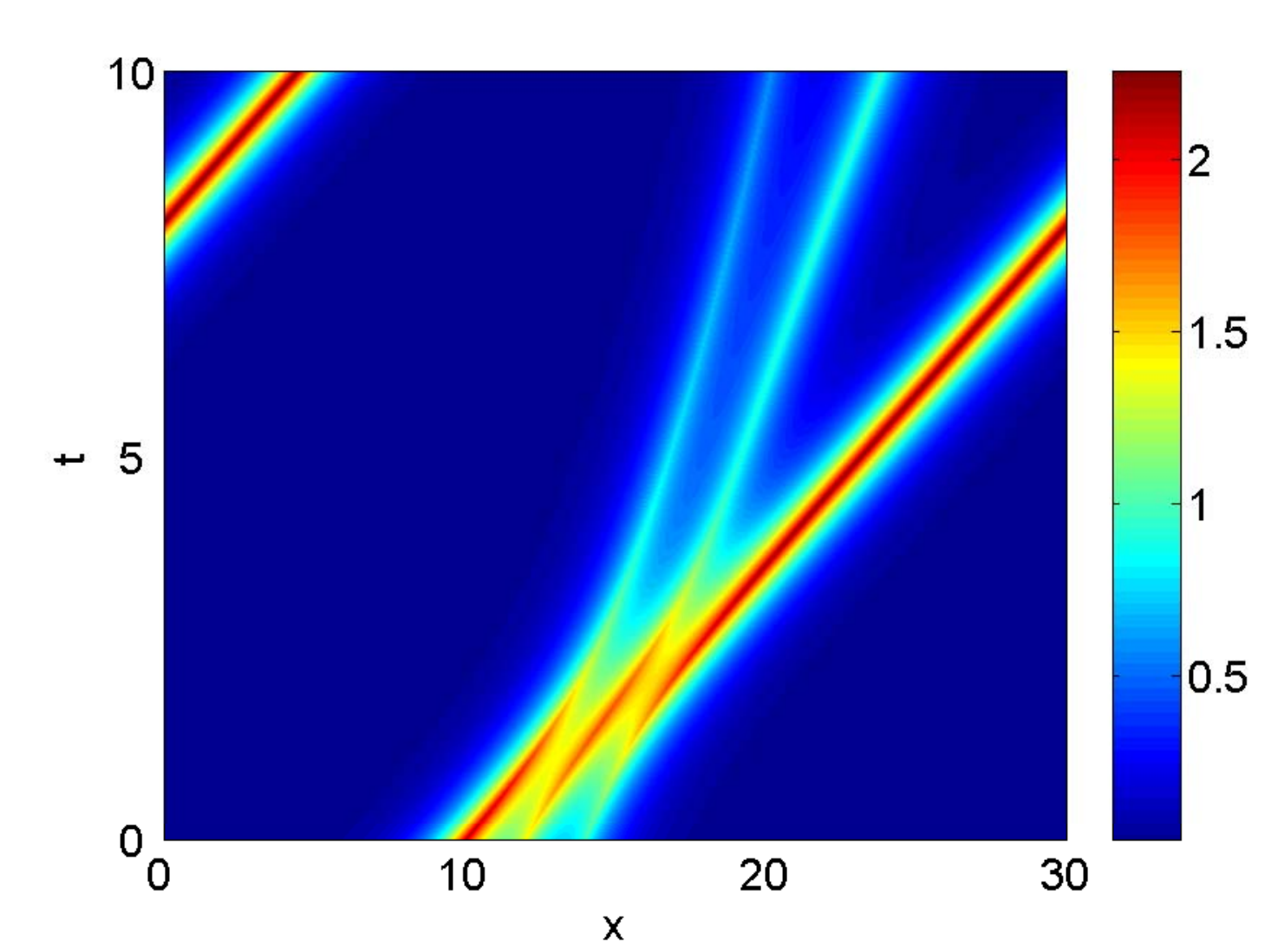}
\end{minipage}
\caption{The three-peakon interaction of the CH equation \eqref{ch:eq:1.1} under $N$=2048 and $\tau$=0.0001.}\label{Fig:ch:5}
\end{figure}

Subsequently, we consider the three-peakon interaction of the CH equation \eqref{ch:eq:1.1} with the initial condition \cite{XS08}
\begin{align*}
u_0(x)=\phi_1(x)+\phi_2(x)+\phi_3(x),\ 0\le x\le 30,
\end{align*}
where
\begin{align*}
\phi_i(x)=\left\lbrace
  \begin{aligned}
  &\frac{c_i}{\cosh(L/2)}\cosh(x-x_i),\ |x-x_i|\le L/2,\\
  &\frac{c_i}{\cosh(L/2)}\cosh(30-(x-x_i)),\ |x-x_i|> L/2,\
  \end{aligned} i=1,2,3.
  \right.\ \ 
  \end{align*}
The parameters are $c_1=2, c_2=1, c_3=0.8, x_1=-5, x_2=-3, x_3=-1,L=30$ and a periodic boundary condition is chosen.
Fig. \ref{Fig:ch:5} shows the contour plot of three-peakons interaction, which shows that the moving peak interaction is resolved very well. The errors in invariants over the time interval $t\in[0,10]$ are displayed in Fig. \ref{Fig:ch:5}, which demonstrates that our scheme has a good conservation of the invariants.

\begin{figure}[H]
\centering\begin{minipage}[t]{60mm}
\includegraphics[width=65mm]{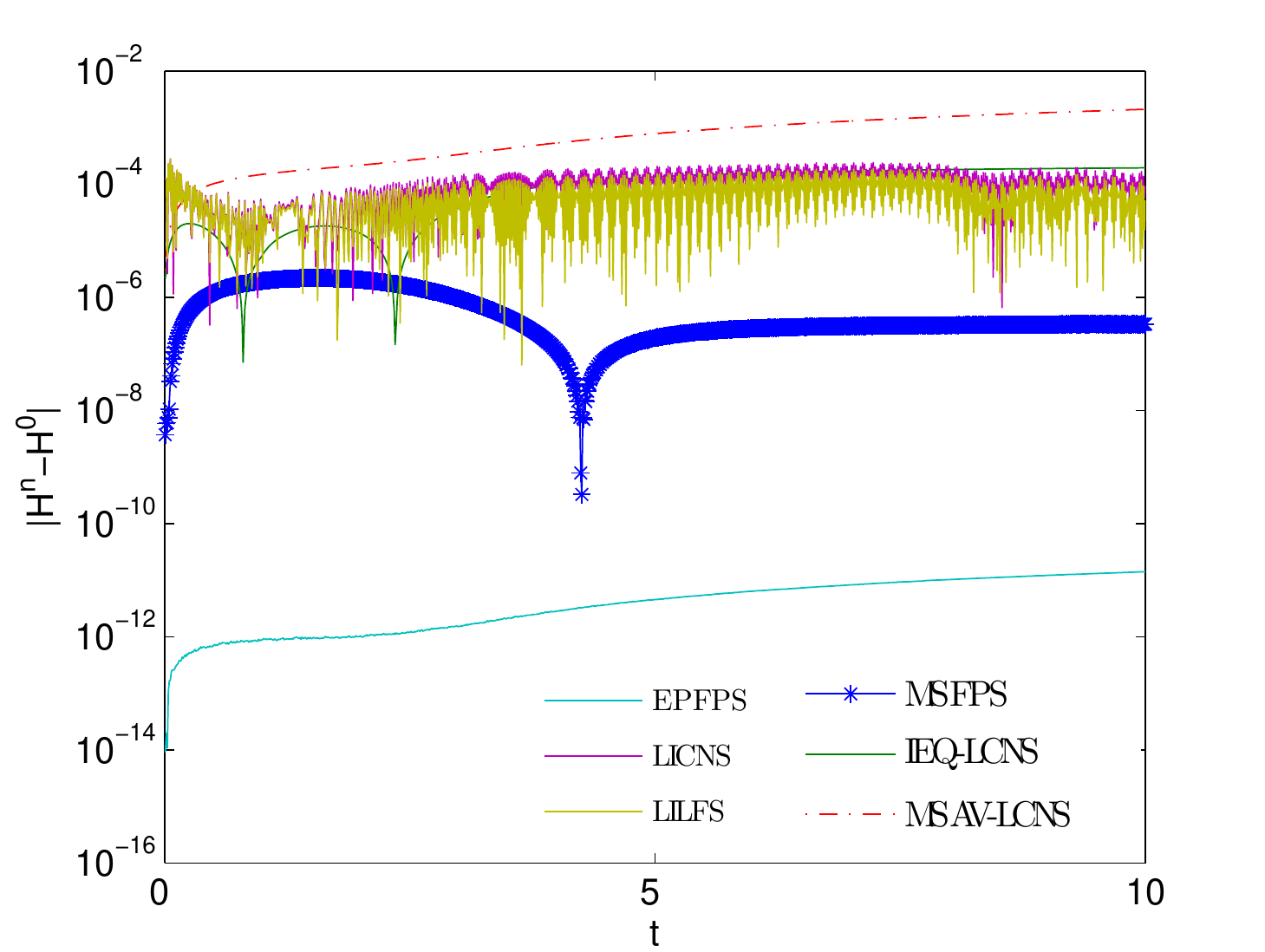}
\caption*{(a) Hamiltonian energy}
\end{minipage}\ \
\begin{minipage}[t]{60mm}
\includegraphics[width=65mm]{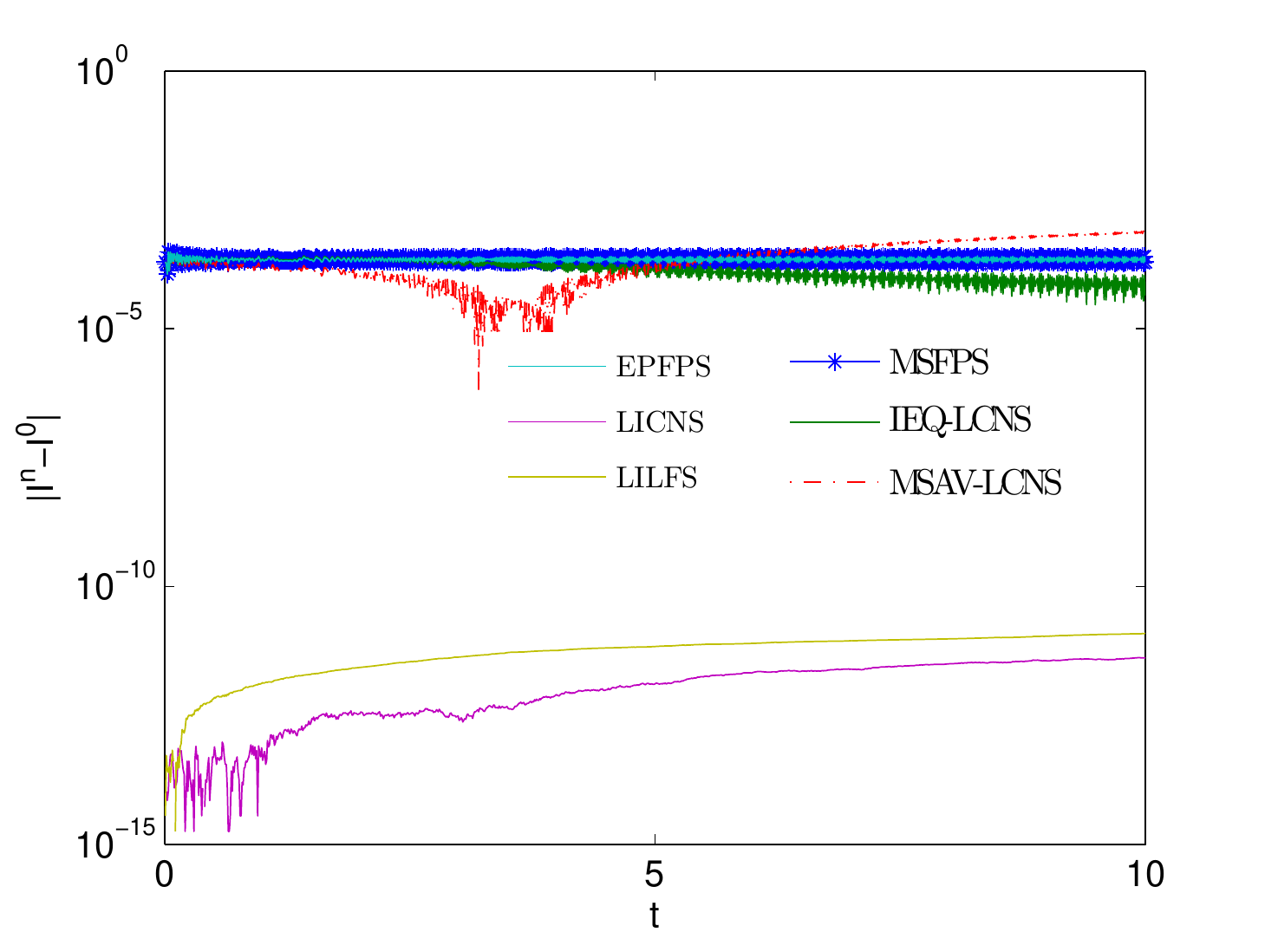}
\caption*{(b) Momentum}
\end{minipage}
\centering\begin{minipage}[t]{60mm}
\includegraphics[width=65mm]{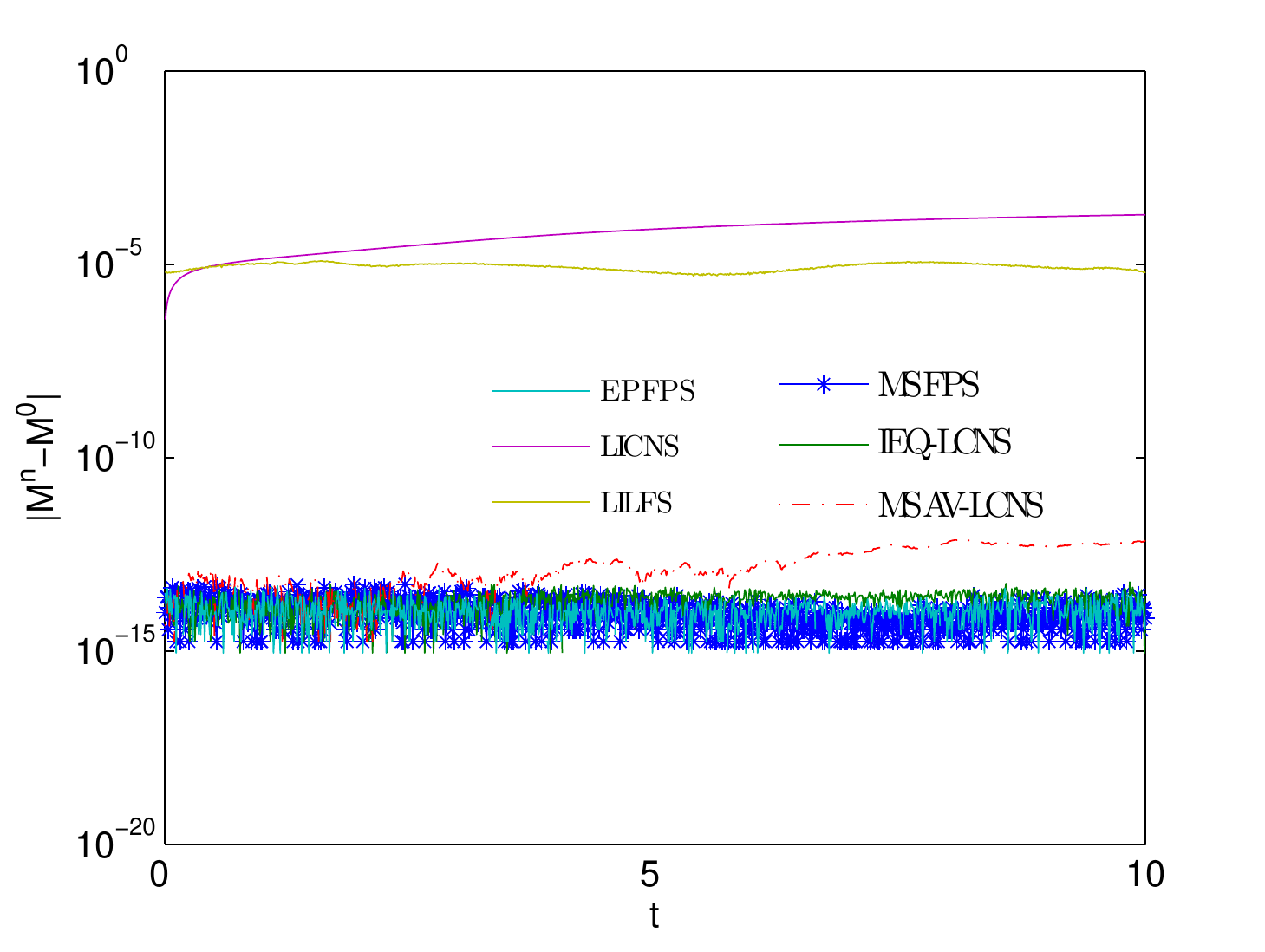}
\caption*{(c) Mass}
\end{minipage}\ \
\begin{minipage}[t]{60mm}
\includegraphics[width=65mm]{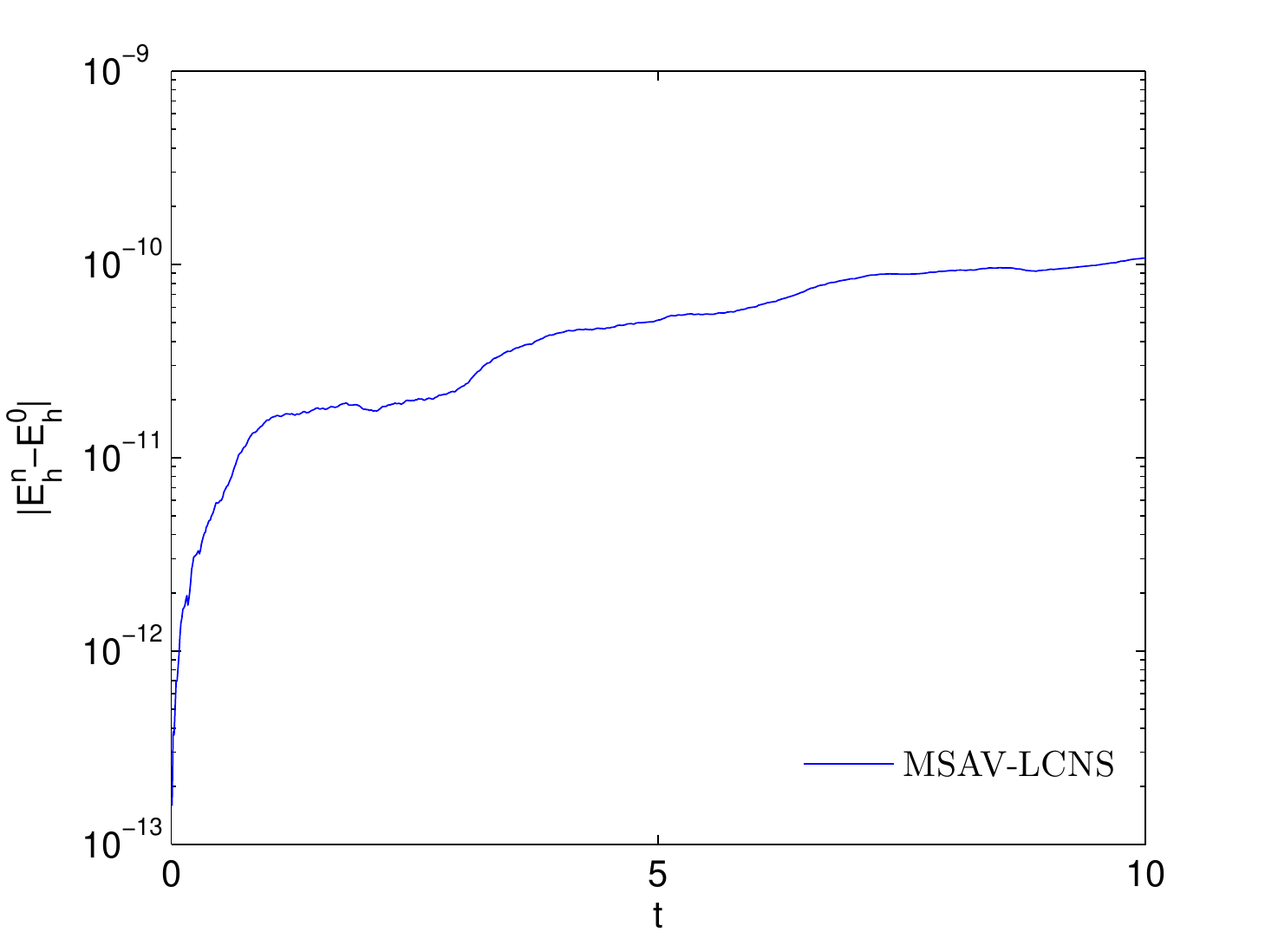}
\caption*{(d) Energy \eqref{ch:eq-4.4}}
\end{minipage}
\caption{The errors in invariants under $N$=2048 and $\tau$=0.0001 over the time interval $t\in[0,10]$.}\label{Fig:ch:6}
\end{figure}

\subsection{A solution with a discontinuous derivative}
Finally, we consider the following initial condition, which has a discontinuous derivative \cite{XS08}
\begin{align*}
u_0(x)=\frac{10}{(3+|x|)^2},\ -30\le x\le 30,
\end{align*}
with a periodic boundary condition.

\begin{figure}[H]
\centering\begin{minipage}[t]{80mm}
\includegraphics[width=80mm]{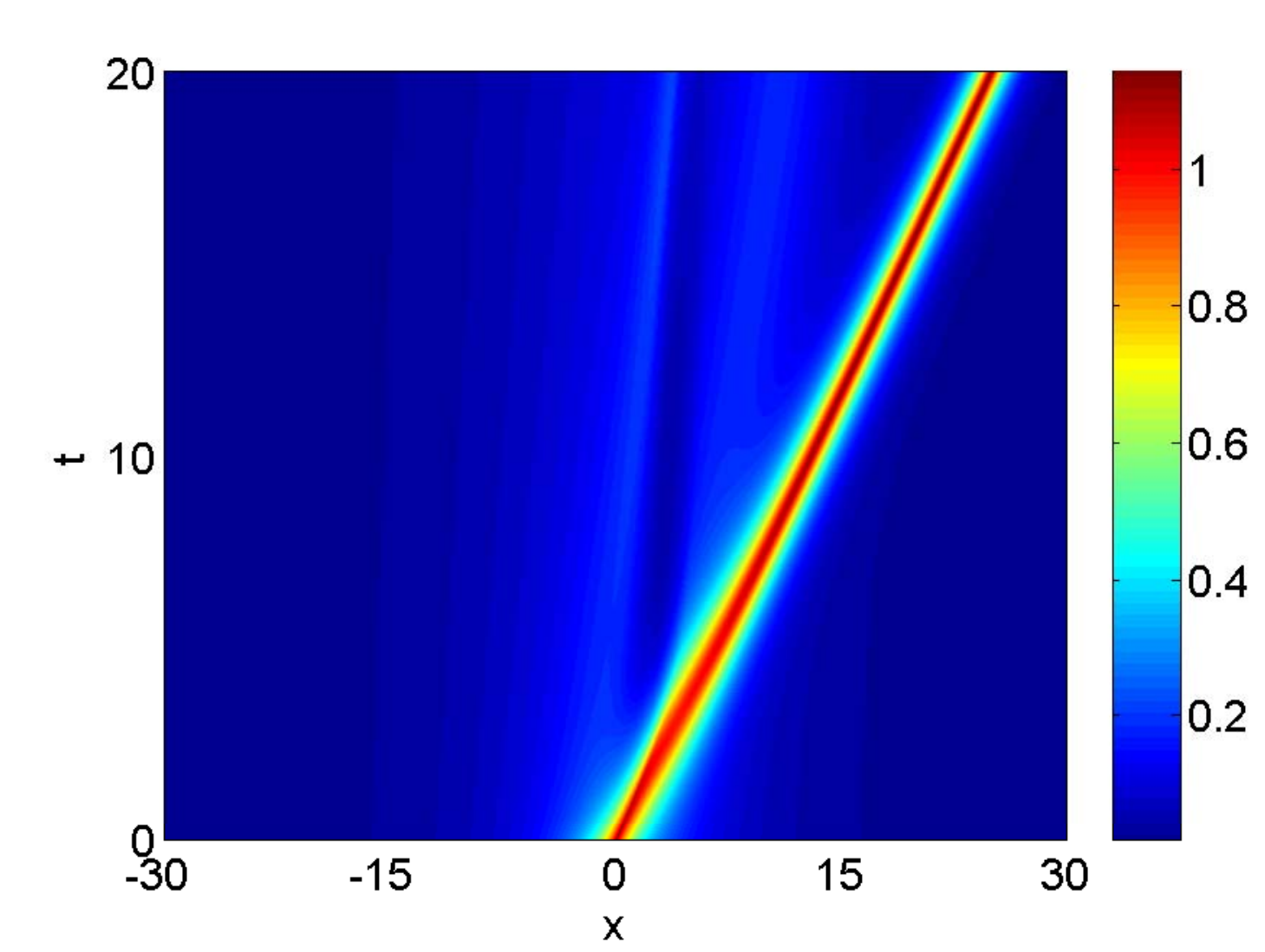}
\end{minipage}\ \
\caption{The solution with discontinuous derivative of the CH equation \eqref{ch:eq:1.1} under $N$=1024 and $\tau$=0.001.}\label{Fig:ch:7}
\end{figure}

\begin{figure}[H]
\centering\begin{minipage}[t]{60mm}
\includegraphics[width=65mm]{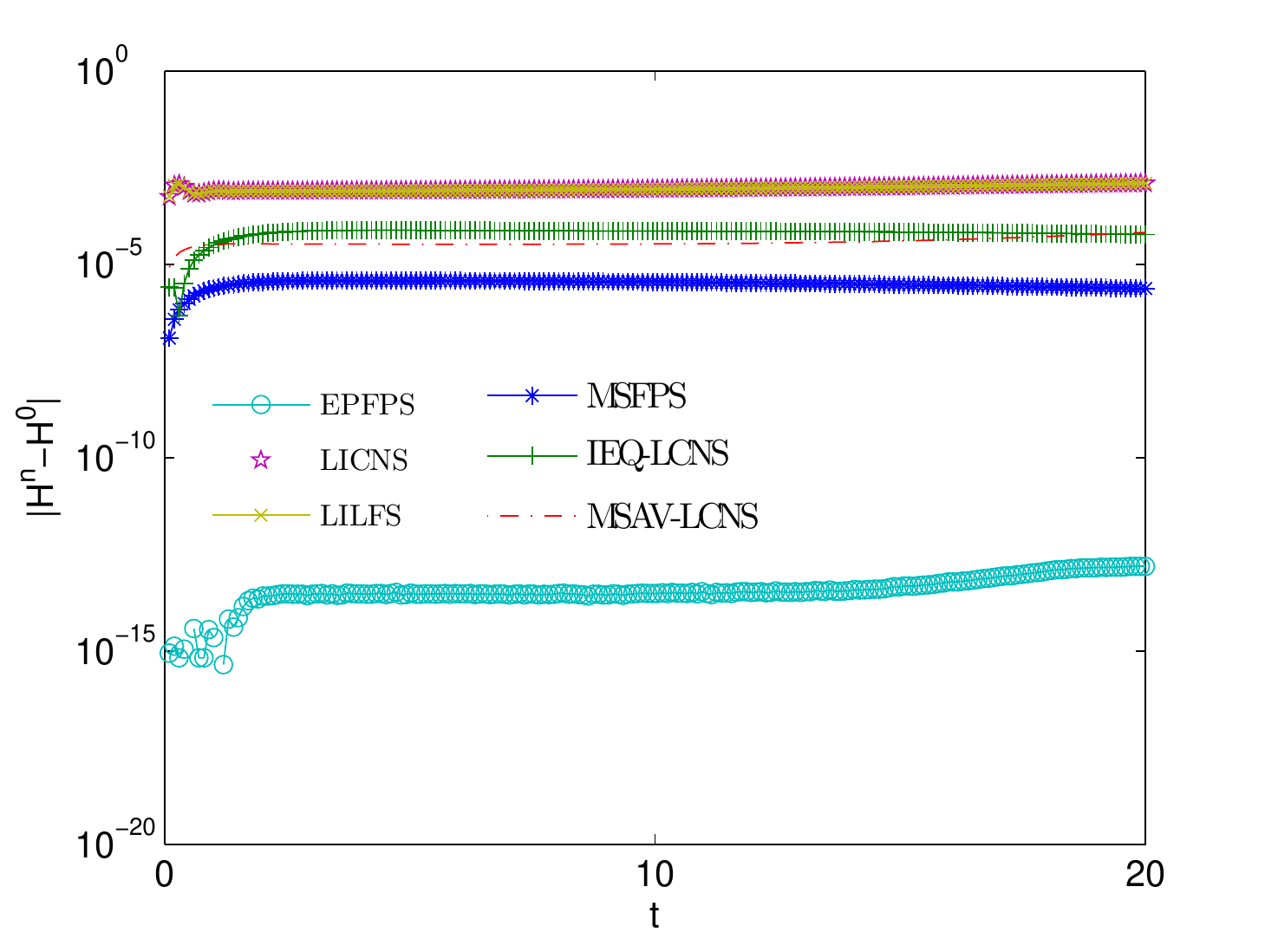}
\caption*{(a) Hamiltonian energy}
\end{minipage}\ \
\begin{minipage}[t]{60mm}
\includegraphics[width=65mm]{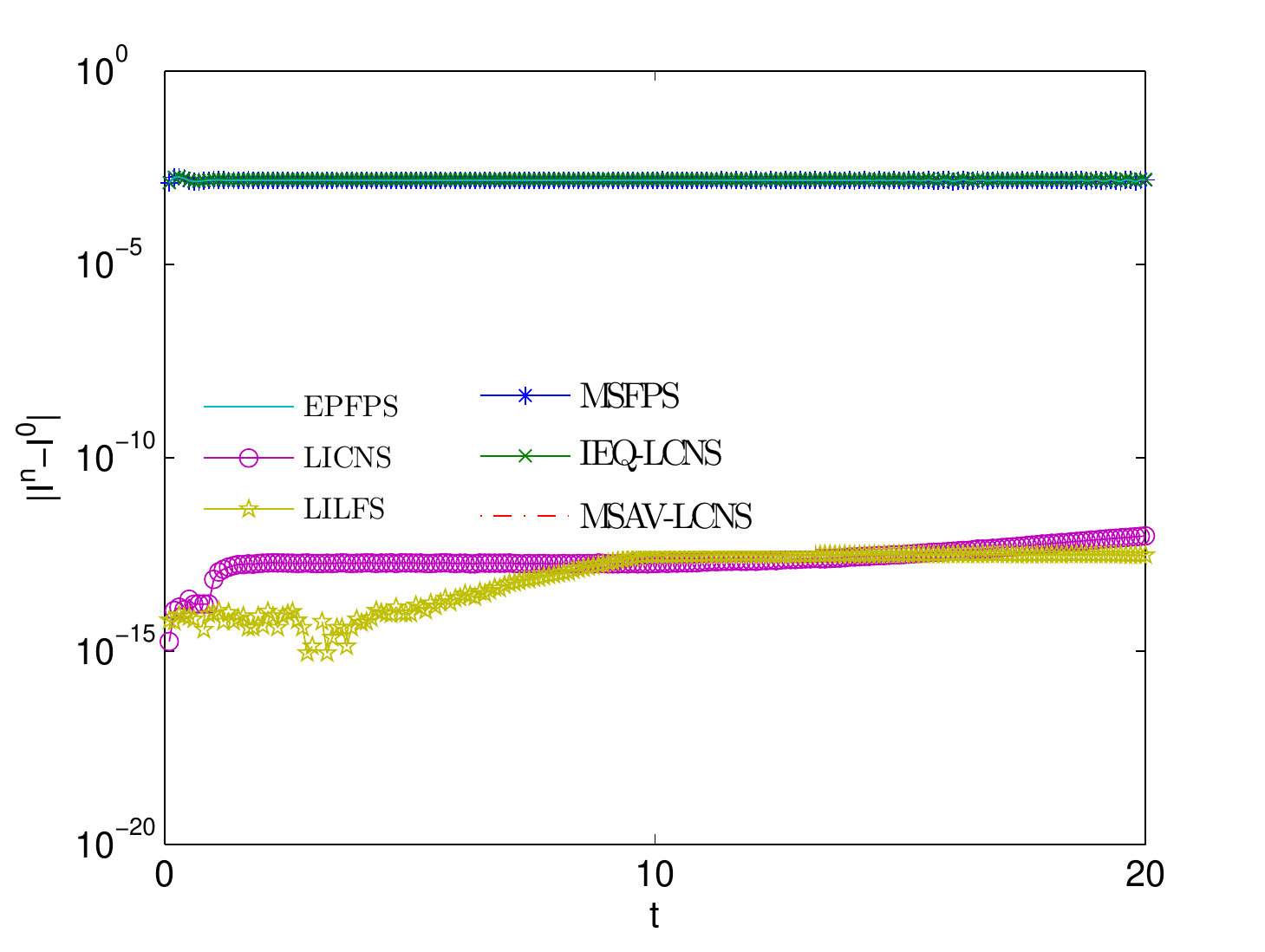}
\caption*{(b) Momentum}
\end{minipage}
\centering\begin{minipage}[t]{60mm}
\includegraphics[width=65mm]{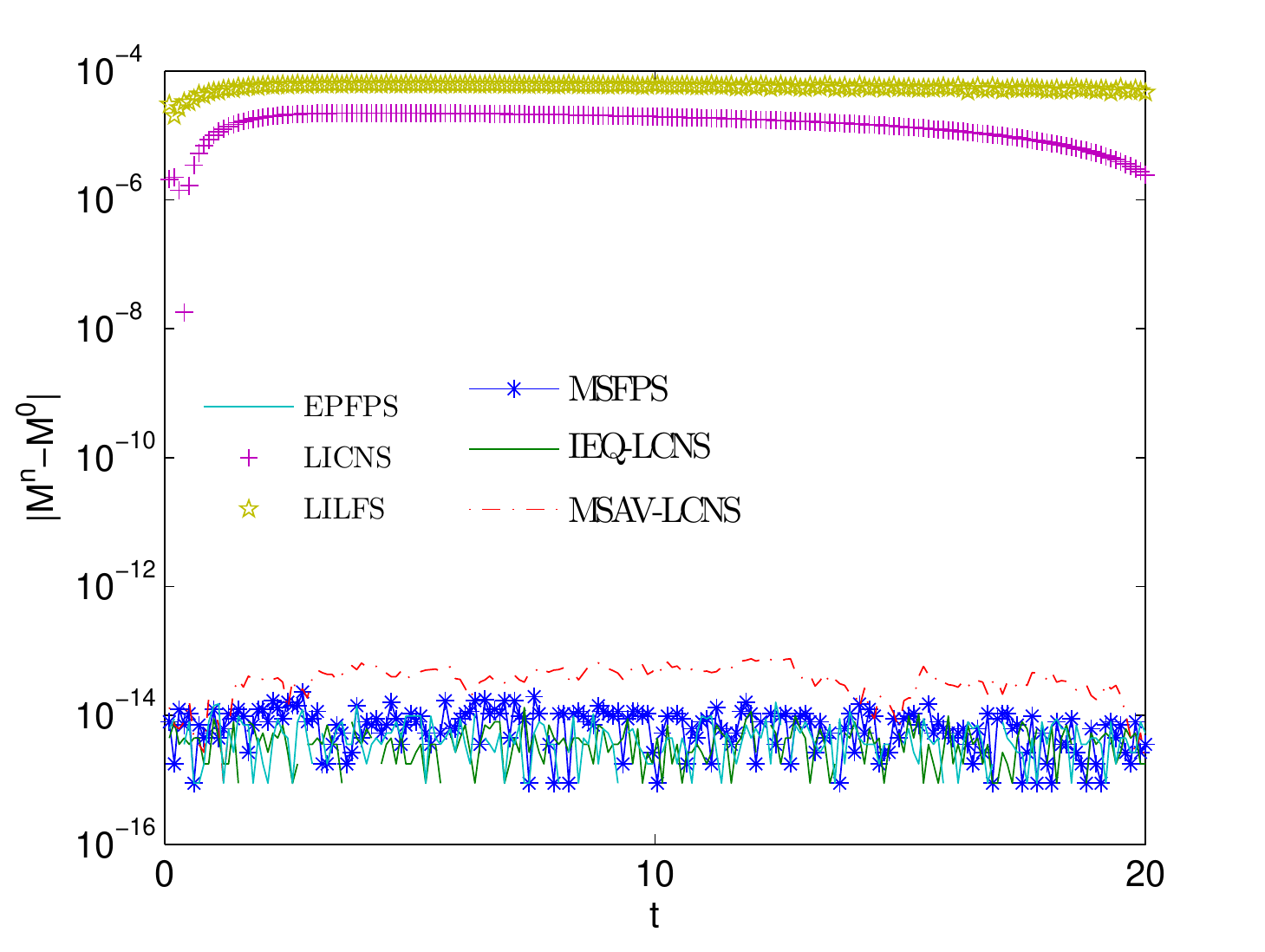}
\caption*{(c) Mass}
\end{minipage}\ \
\begin{minipage}[t]{60mm}
\includegraphics[width=65mm]{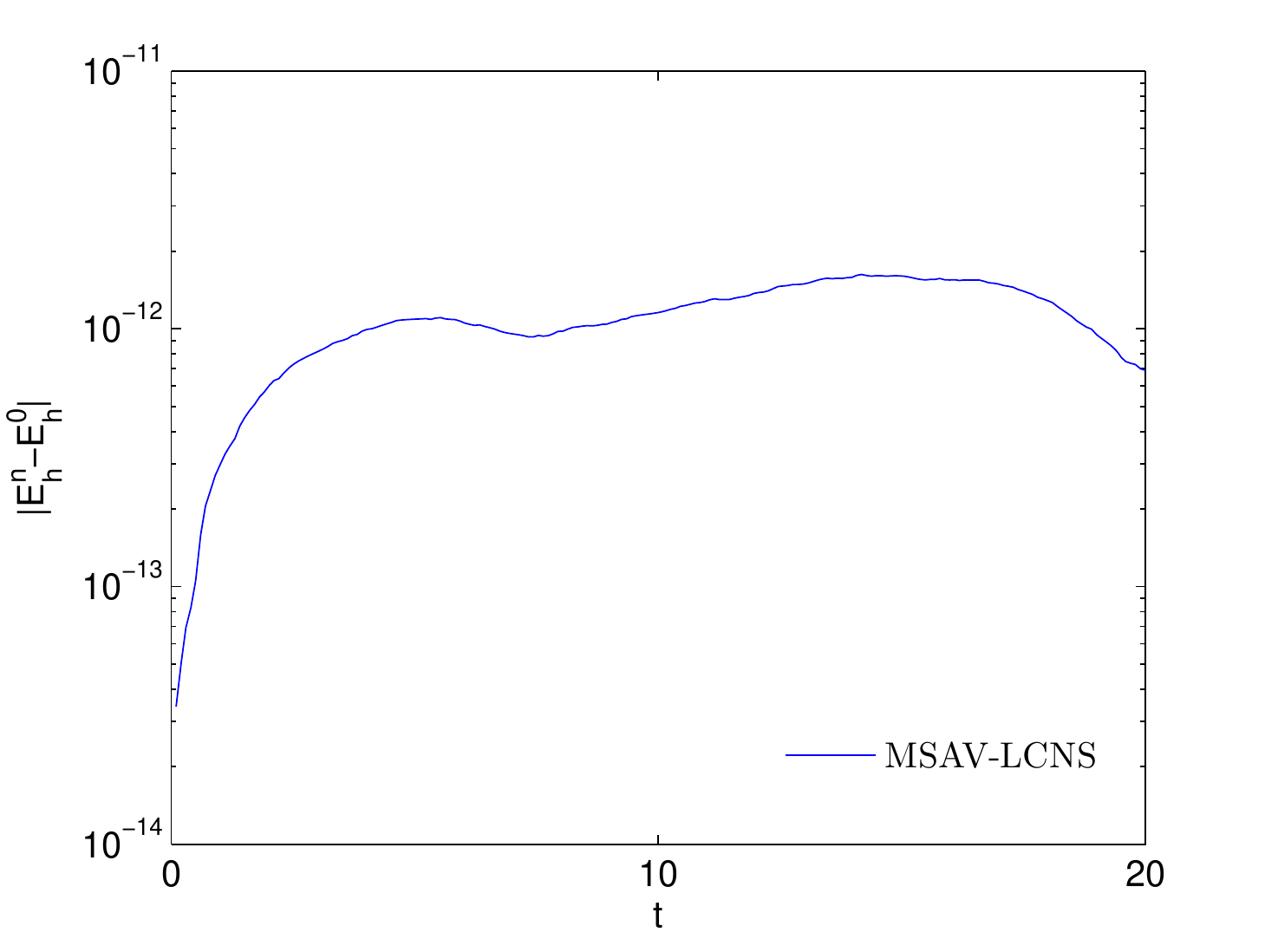}
\caption*{(d) Energy \eqref{ch:eq-4.4}}
\end{minipage}
\caption{The errors in invariants under $N$=1024 and $\tau$=0.001 over the time interval $t\in[0,20]$.}\label{Fig:ch:8}
\end{figure}

Fig. \ref{Fig:ch:7} shows the contour plot of the solutions with discontinuous derivative. Fig. \ref{Fig:ch:8} shows the errors in invariants over the time interval $t\in[0,20]$. From Figs. \ref{Fig:ch:7} and \ref{Fig:ch:8}, it is clearly demonstrated that the proposed scheme has a good resolution of the solution comparable with that in Refs. \cite{HR06,XS08}, and can preserve the modified energy exactly.

\section{Concluding remarks}\label{Sec:Ch:6}
In this paper, we present a novel linearization (energy
quadratization) strategy to develop second order, fully discrete,
linearly implicit scheme for the CH equation \eqref{ch:eq:1.1}. The proposed scheme is proven to preserve the discrete modified energy and enjoys the same computational advantages as the schemes obtained by the classical SAV approach. Several numerical examples are presented to illustrate
the efficiency of our numerical scheme. Comparing with some existing structure-preserving schemes of same order in both time and space, our scheme shows remarkable efficiency. The presented strategy can be directly extended to propose linearly implicit energy-preserving schemes for a broad class of energy-conserving systems, such as the KdV equation, etc. However, to the best of our knowledge, the construction of arbitrarily high-order linearly implicit energy-preserving schemes is still not available for the CH equation \eqref{ch:eq:1.1}, which is an interesting topic for future studies.

\section*{Acknowledgments}

The authors would like to express sincere gratitude to the referees for their insightful
comments and suggestions. Chaolong Jiang's work is partially supported by the National Natural Science Foundation of China (Grant No. 11901513), the Yunnan Provincial Department of Education Science Research Fund Project (Grant No. 2019J0956) and the Science and Technology Innovation Team on Applied Mathematics in Universities of Yunnan. Yuezheng Gong's work is partially supported by the Natural Science Foundation of Jiangsu Province
(Grant No. BK20180413) and the National Natural Science Foundation of China (Grant No. 11801269). Wenjun Cai's work is partially supported
by the National Natural Science Foundation of China (Grant No. 11971242) and the National
Key Research and Development Project of China (Grant Nos. 2018YFC0603500, 2018YFC1504205).
 Yushun Wang's work is partially supported by the National Natural Science Foundation of China (Grant No. 11771213).


\bibliographystyle{plain}

\end{document}